\documentclass[11pt,times]{article}

\usepackage{graphicx,color}

\usepackage[colorlinks=true, breaklinks=true, pdfstartview=FitV, linkcolor=red, citecolor=blue, urlcolor=blue]
{hyperref}

\newcommand{\Vec}[1]{\mbox{\boldmath$#1$}}

\begin{document}

\title{Nonlinear librations of distant retrograde orbits: a perturbative approach --The Hill problem case\thanks{%
Nonlinear Dyn (2018). \href{https://doi.org/10.1007/s11071-018-4304-0}{https://doi.org/10.1007/s11071-018-4304-0} (Pre-print version)}
}

\author{Martin Lara\thanks{%
GRUCACI, University of La Rioja, 26006 Logro\~no, La Rioja, Spain. {\tt mlara0@gmail.com}
}
}

\date{\color{red}draft of January 8, 2018}

\maketitle

\begin{abstract}
The non-integrability of the Hill problem makes that its global dynamics must be necessarily approached numerically. However, the analytical approach is feasible in the computation of relevant solutions. In particular, the nonlinear dynamics of the Hill problem close to the origin, and the libration point dynamics have been thoroughly investigated by perturbation methods. Out of the Hill sphere, the analytical approach is also feasible, at least in the case of distant retrograde orbits. Previous analytical investigations of this last case succeeded in the qualitative description of the dynamics, but they commonly failed in providing accurate results. This is a consequence of the essential dependance of the dynamics on elliptic functions, a fact that makes to progress in the perturbation approach beyond the lower orders of the solution really difficult. We propose an alternative perturbation approach that allows us to provide a very simple low order analytical solution in trigonometric functions, on the one hand, and, while still depending on special functions, to compute higher orders of the solution, on the other.
\end{abstract}

\section{Introduction}

The Hill problem is a limit case of the three-body problem in which two of the involved masses are very small when compared to the third, dominant mass \cite{HenonPetit1986}. If, besides, one of the small masses is negligible and moves close enough to the other small mass, the Hill problem can be viewed as a particular case of the restricted three body problem \cite{Szebehely1967}. The fact that, after a convenient nondimensionalization of the physical units, the Hill problem does not depend on any physical parameter endows this dynamical model with a wide generality. Indeed, the Hill problem is representative of the dynamics of different astronomical and astrodynamical problems, which range from the Sun-Earth-Moon motion, Hill's original application \cite{Hill1878}, to the leader-follower motion of artificial satellites \cite{ClohessyWiltshire1960}. The range of applications of the Hill problem include also the study of dust interaction in planetary rings \cite{PetitHenon1986}, the description of the motion around planetary satellites \cite{LaraRusselVillacl2007Meccanica}, or different aerospace engineering applications, which, in particular, include those related to the the libration point dynamics \cite{Gomezetal2003}, to just mention a few.
\par

When the small masses are far enough away from each other, the Hill problem can be simplified to the relative motion of two particles each one moving in a Keplerian orbit about the bigger mass, which, therefore, becomes an integrable problem. However, when the mutual attraction of the two small masses is comparable to the differential attraction of the dominant mass the Hill problem turns into non-integrable. The non-integrable problem admits two symmetric equilibria, the so-called libration points, as well as a variety of periodic orbit families. Further than these particular solutions, the global dynamics of the Hill problem must unavoidably be approached with numerical techniques \cite{Henon1969,SimoStuchi2000,Zotos2017}.
\par

On the other hand, analytical approximations to the Hill problem have been computed in particular cases of relevant interest. Indeed, the Hamiltonian of the Hill problem can be rearranged as a perturbed Keplerian problem when the motion of the particle with negligible mass occurs inside the gravitational sphere of influence, or Hill's sphere, of the other small body \cite{LaraPalacianRussell2010}. This rearrangement makes the inclusion of additional effects in the perturbation model natural, in order to provide better approximations to the real dynamics \cite{LidovYarskaya1974,LaraSanJuanFerrer2005,SanJuanLaraFerrer2006}.
\par

Another case that admits the perturbation approach appears in the study of the dynamics about the libration points. Indeed, after translation of the origin to the libration point, followed by the standard expansion of the Keplerian potential in Legendre polynomials, the Hill problem Hamiltonian takes the form of a hyperbolic term plus a perturbed elliptic oscillator. Further than the traditional reduction to the center manifold that removes the hyperbolic components of the motion \cite{GomezMarcoteMondelo2005}, or the Lindstedt-Poincar\'e computation of particular periodic orbits \cite{ZagourasMarkellos1985}, the relevant dynamics in this region of the Hill problem is captured with a single normal form Hamiltonian. This reduced Hamiltonian not only discloses all the major features of the qualitative dynamics about the libration points \cite{Lara2017}, but it also provides the higher order terms that are required in the accurate computation of the more relevant solutions \cite{LaraPerezLopez2017}. 
\par

A third instance in which the analytical approach is feasible is the case of co-orbital motion of the two small primaries \cite{Benest1976,Namouni1999}. Now, the Keplerian potential is taken as a perturbation of the (integrable) quadratic Hamiltonian representing the relative motion of the smaller bodies, each of which evolves with Keplerian motion about the massive body. Then, the linearized dynamics of the Hill problem reduces to a moving ellipse whose major axis doubles the length of the minor axis, and whose center evolves with linear motion in the $y$ axis direction ---along which the major axis lies. Nonlinearities introduced by the gravitation of the small, non-vanishing mass, may modify this linear displacement to the extent of turning the motion of the center of the reference ellipse into long-period librations.
\par

In this last case, the perturbation solution depends inherently on elliptic functions, a fact that  complicates considerably extending the computations to the higher orders that may be required in practice \cite{LidovVashkovyak1993,LidovVashkovyak1994}. Furthermore, the first order approach is commonly limited to the computation of the secular terms of the solution. In consequence, the analytical approach is constrained to the qualitative description of the motion. Due to this lack of precision of current analytical solutions, numerical explorations based on the computation of periodic orbits as well as other invariant manifolds of restricted models, or the use of Lyapunov indicators, are usual approaches in the study of the so-called quasi-satellite, or distant retrograde orbits \cite{LamWhiffen2005,LaraRussellVillac2007,LaraRusselVillacl2007Meccanica,MingShijie2009,GilSchwartz2010,StramacchiaColomboBernelliZazzera2016}. Still, higher order series expansions are available, but restricted to some particular situations \cite{HenonPetit1986}.
\par

Recent discoveries of different co-orbital asteroids in quasi-satellite orbits (see \cite{delaFuenteMarcosCyR2014} and references therein) as well as proposed missions requiring  long duration quarantine orbit, such as the Asteroid Redirect Mission \cite{Abelletal2017} or the recent DePhine proposal to the European Sapce Agency \cite{OberstWillnerWickhusen2017}, stir new interest in the study of distant retrograde orbits \cite{Sidorenkoetal2014,BezroukParker2017,Perozzietal2017,PousseRobutelVienne2017}. This fact motivates us for further digging into the possibilities offered by the analytical approach.
\par

We revisit the distant retrograde orbits problem in the approximation to the dynamics provided by the planar case of the Hill problem, and propose an alternative perturbation solution that allows the computation of higher orders of the reduced Hamiltonian beyond the usual first order approach. For the lower orders of the solution, the Hamiltonian flow can be integrated analytically leading to an explicit solution that is free from special functions. The accuracy of this lower order solution is comparable to previous solutions in the literature, yet it provides a deeper insight on the dynamics as well as easier evaluation.
\par

The lower order solution is useful for capturing the periodic or quasi-periodic motion which is typical of distant retrograde orbits, but it fails in the case of higher librations of the reference ellipse, which may cause close encounters between the small masses. When the Hamiltonian reduction is extended to higher orders, corresponding Hamilton equations are still solved by quadrature. However, the closed form of these quadratures, if found, would require the use of special functions in an intricate way. Besides, the following inversion of the solution would be required to make it explicit. These facts not only deprive the analytical solution of the expected insight, but they also make its evaluation difficult. Still, useful approximations to the higher order analytical solution can be computed using Lindstedt series (see \cite{Murdock1999}, for instance). However, the increasing complexity of the perturbation approach when reaching higher orders, which requires dealing with the analytical integration of special functions, establishes a practical limit in the analytic computations which prevents the correct modeling of the dynamics of orbits with large librations. Alternatively, the numerical integration of the Hamilton equations of the secular Hamiltonian is carried out with very large step sizes, and hence is very fast and efficient. With this last approach, additional higher order effects of the gravitational attraction of the smaller primary can be incorporated, on average, to the perturbation solution.
\par

The paper is organized as follows. The Hamiltonian formulation of the Hill problem, which is rearranged as a perturbation problem valid for the case of distant retrograde orbits, is presented in Section \ref{s:problem}, in which the linearized dynamics is discussed after completely reducing the quadratic Hamiltonian, and the nonlinear, perturbation term is expanded in the canonical variables that reduce the unperturbed problem. The low order, explicit solution in trigonometric functions is derived in Section \ref{s:low}. Finally, Section \ref{s:high} explores the validity of higher orders of the perturbation solution in the computation of quasi-satellite orbits with large librations.

\section{Hill's equations for relative motion} \label{s:problem}

In spite of its wider generality \cite{HenonPetit1986}, the Hill problem is commonly presented as a particular case of the restricted three-body problem. In that context, the Hill problem describes the motion of a massless body under the gravitational influence of two massive bodies, called the primaries, under the assumptions that: 1) the mass of one of the primaries is much bigger than the mass of the other, and 2) the distance between the primaries is considerably larger than the distance of the massless body to the origin, which is taken in the smaller primary. Under the additional assumption of circular motion of the primaries, the equations of motion of the massless body can be derived from the Hamiltonian of the Hill problem \cite{Szebehely1967,BoccalettiPucacco1998v1}
\begin{equation}
\mathcal{H}=\frac{1}2\,(\Vec{X}\cdot\Vec{X})-\Vec\omega\cdot(\Vec{x}\times\Vec{X})+\mathcal{V}(\Vec{x};\Vec{\omega}),
\end{equation}
where $\Vec{x}$ is the position vector of the massless body in a rotating frame centered at the smaller primary, with the $x$ axis in the direction from the bigger to the smaller primary, the $z$ axis is defined by the direction of the rotation of the system $\Vec{\omega}$, and the $y$ axis completes a direct orthogonal frame; $\Vec{X}$ notes the conjugate momentum to $\Vec{x}$, and the potential $\mathcal{V}$ is
\begin{equation}
\mathcal{V}=-\frac{\mu}r+\frac{1}2\omega^2r^2\left(1-3\frac{x^2}{r^2}\right),
\end{equation}
where $r=\|\Vec{x}\|$, $\omega=\|\Vec\omega\|$, and $\mu$ is the gravitational parameter of the smaller primary.
\par

Units of length and time are customarily chosen in such a way that $\mu$ and $\omega$ become the unity, showing that the Hill problem does not depend on physical parameters. However, we maintain both parameters in following derivations to avail ourselves with an additional test on the correctness of the analytical developments by checking dimensions.
\par

Besides, we constrain ourselves to the case o planar motion in the $x,y$ plane, in which the Hamiltonian of the Hill problem is written in scalar variables as
\begin{equation} \label{planarHam}
\mathcal{H}=\frac{1}2\left(X+\omega{y}\right)^2+\frac{1}2(Y-\omega{x})^2-\frac{3}{2}\omega^2x^2-\frac{\mu}r.
\end{equation}
\par

Furthermore, we are only interested in the particular case of co-orbital motion out of the Hill sphere of influence of the smaller primary. In that case, the gravitational attraction of the smaller primary, represented by the Keplerian potential $-\mu/r$, will be small. Therefore, Eq.~\ref{planarHam} can be set in the perturbative arrangement
\begin{equation}
\mathcal{H}=\mathcal{H}_0+\epsilon\mathcal{H}_1,
\end{equation}
where
\begin{equation} \label{H0}
\mathcal{H}_0=\frac{1}2\left(X+\omega{y}\right)^2+\frac{1}2(Y-\omega{x})^2-\frac{3}{2}\omega^2x^2,
\end{equation}
is integrable, 
\begin{equation} \label{H1}
\mathcal{H}_1=-\frac{\mu}r,
\end{equation}
and $\epsilon$ is a formal small parameter which is used to manifest that the effect of 
$\mathcal{H}_1$, the perturbation, is much smaller than $\mathcal{H}_0$.

\subsection{Linearized dynamics: complete reduction}

In preparation of the perturbation approach, the integration of the linearized dynamics derived from the quadratic Hamiltonian (\ref{H0}) is carried out by complete reduction \cite{Arnold1989}.
\par

Thus, the canonical transformation 
\[
\mathcal{T}:(x,y,X,Y)\longrightarrow(\phi,q,\Phi,Q;\omega),
\]
given by the equations
\begin{equation} \label{tdir}
\begin{array}{rcl}
\omega{x} &=&\displaystyle 2Q+\sqrt{2\omega\Phi}\sin\phi, \\[1ex]
\omega{y} &=&\displaystyle \omega{q}+2\sqrt{2\omega\Phi}\cos\phi, \\[1ex]
 X &=&\displaystyle -\omega{q}-\sqrt{2\omega\Phi}\cos\phi, \\[1ex]
 Y &=&\displaystyle -Q-\sqrt{2\omega\Phi}\sin\phi,
\end{array}
\end{equation}
converts Eq.~(\ref{H0}) into the completely reduced Hamiltonian in the new variables 
\begin{equation} \label{quadratic}
\mathcal{H}_0\circ\mathcal{T}\equiv\mathcal{K}_0=\omega\Phi-\mbox{$\frac{3}{2}$}Q^2,
\end{equation}
where the angle $\phi$ and the distance $q$ are ignorable coordinates. In consequence, the action $\Phi$ and the velocity $Q$ are constant, whereas $\phi$ and $q$ evolve linearly with time, viz.
\begin{eqnarray}
\dot\phi &=& \frac{\partial\mathcal{H}_0}{\partial\Phi}= \omega, \\
\dot{q} &=& \frac{\partial\mathcal{H}_0}{\partial{Q}}=-3Q,
\end{eqnarray}
as derived from Hamilton equations.
\par

Therefore, the Hamiltonian flow derived from Eq.~(\ref{quadratic}) is trivially integrated, to give
\begin{eqnarray} \label{solphi}
\phi &=& \phi_*+\omega{t}, \\  \label{solq} 
q &=& q_*-3Qt,
\end{eqnarray}
in which the subindex $*$ is used to note the value of a variable at the time $t=0$.
\par

The solution is obtained in the original variables after replacing Eqs.~(\ref{solphi}) and (\ref{solq}) into the direct transformation given in Eq.~(\ref{tdir}). It yields
\begin{eqnarray} \label{qs0x}
x &=& 2Q/\omega+K_1\cos\omega{t}+K_2\sin\omega{t}, \\ \label{qs0y}
y &=& q_*-3Qt+2K_2\cos\omega{t}-2K_1\sin\omega{t},
\end{eqnarray}
as well as corresponding expressions for their conjugate momenta $X$ and $Y$, where 
\begin{eqnarray*}
K_1 &=& \sqrt{2\Phi/\omega}\sin\phi_*, \\
K_2 &=& \sqrt{2\Phi/\omega}\cos\phi_*,
\end{eqnarray*}
and the arbitrary constants $\Phi$, $Q$, $\phi_*$, and $q_*$ are obtained as functions of the initial conditions using the inverse transformation of Eq.~(\ref{tdir}), which is easily obtained from the obvious relations
\begin{equation} \label{tinv}
\begin{array}{rcl}
\omega{x}+Y &=& Q, \\[1ex]
\omega{y}+X &=& \sqrt{2\omega\Phi}\cos\phi, \\[1ex]
\omega{x}+2Y &=& -\sqrt{2\omega\Phi}\sin\phi, \\[1ex]
\omega{y}+2X &=& -\omega{q}.
\end{array}
\end{equation}
We obtain,
\begin{equation}\label{iicc} \begin{array}{rcl}
Q &=& Y_*+\omega{x}_*, \\
K_1 &=& -x_*-2Y_*/\omega, \\
K_2 &=& y_*+X_*/\omega, \\
q_* &=& -y_*-2X_*/\omega.
\end{array}
\end{equation}
\par

As expected, the solution in Eqs.~(\ref{qs0x})--(\ref{qs0y}) is the same solution previously given by Benest \cite{Benest1976}, who based his computations in the variation of parameters approach. Namely, in configuration space, the orbit is a moving ellipse with semi-major axis of length
\begin{equation} \label{semimajor}
A=2\sqrt{2\Phi/\omega},
\end{equation}
lying along the $y$ axis direction, and semi-minor axis of length
\begin{equation} \label{semiminor}
B=\sqrt{2\Phi/\omega}=A/2.
\end{equation}
The instantaneous center of the ellipse has coordinates
\begin{eqnarray}
x_\mathrm{C} &=& 2Q/\omega, \\
y_\mathrm{C} &=& q,
\end{eqnarray}
where the abscissa is constant and the ordinate varies linearly with time, as given by Eq.~(\ref{solq}), while $\phi$ gives the phase of the ellipse in non-rotating coordinates, cf.~\cite{Benest1976}.
\par

\begin{figure}[htb]
\begin{center}
\includegraphics[scale=0.75, angle=0]{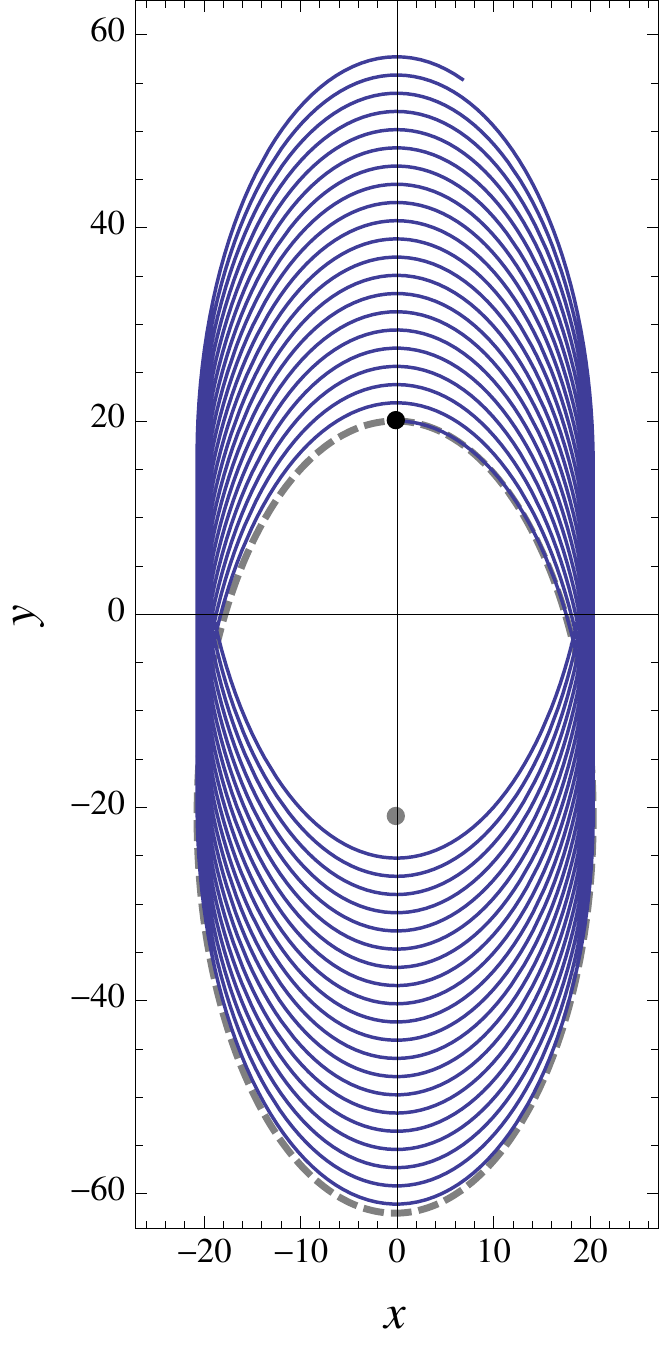}
\caption{Sample quasi-satellite orbit (blue line) with the initial reference ellipse (gray, dashed) superimposed. The black dot marks the initial point and the gray one the center of the initial ellipse.}
\label{f:QSorbit}
\end{center}
\end{figure}

An example of a quasi-satellite orbit is illustrated in Fig.~\ref{f:QSorbit}, which has been computed from Eqs.~(\ref{qs0x})--(\ref{qs0y}) with $\omega=1$, and has been propagated up to $t=126$ starting from the initial conditions $x_*=0$, $y_*=20$, $X_*=0.5$, $Y_*=-0.1$. The initial reference ellipse (dashed, gray line of the figure, whose center is represented by a gray dot with coordinates $x=0.2$, $y=-21$), has been superimposed to the quasi-satellite orbit. Note that, from the first of Eq.~(\ref{iicc}), $|Q|=|Y_*|<\omega/3$ in this example, and hence the linear displacement of the reference ellipse along the $y$ axis progresses slowly when compared with the fast motion, with frequency $\omega$, of the particle of negligible mass along this moving ellipse.
\par

When initial conditions are chosen such that
\begin{equation} \label{periodic0}
Y_*=-\omega{x}_*,
\end{equation}
then $Q$ and $x_\mathrm{C}$ vanish, and the secular component of the Cartesian coordinate $y$ is removed from Eq.~(\ref{qs0y}). Therefore, the massless body enjoys periodic motion. Besides, from Eq.~(\ref{solq}), the last of Eq.~(\ref{tinv}) is rewritten
\[
\omega{y}+2X=-\omega(q_*-3Qt).
\]
Therefore, if, in addition to Eq.~(\ref{periodic0}), the remaining initial conditions are chosen such that 
\begin{equation} \label{centered0}
X_*=-\frac{1}{2}\omega{y}_*,
\end{equation}
then it happens that $q_*=0$, and the center of the periodic ellipse is located at the origin. This is the typical case of distant retrograde, periodic orbits \cite{Henon1969}.
\par

\subsection{Perturbation arrangement}

When the gravitational attraction of the primary becomes non-negligible, the solution in Eqs.~(\ref{qs0x}) and (\ref{qs0y}) no longer applies, and the nonlinearities of the motion produced by the gravitation of the smaller primary, as given by the Hamiltonian term $\mathcal{H}_1$ in Eq.~(\ref{H1}), must be taken into account. Passing orbits, which arbitrarily depart from the smaller primary, may survive to the perturbation. However, depending on the initial state, the accumulation with time of the small effects produced by the perturbation can be enough to turning the linear displacement of the center of the reference ellipse into a long-period libration about the origin. Therefore, some of the quasi-satellite orbits of the Hill problem will remain trapped in the vicinity of the smaller primary without restricting to periodic motion \cite{Henon1970}.
\par

The reduction of the problem by integrals is not possible now. Still, the use perturbation methods make the computation of useful analytical solutions possible, the validity of which is constrained to some regions of phase space and apply only for a limited time interval ---which can be really long depending on the truncation order of the perturbation solution \cite{Nayfeh2004}. 
\par

Application of the canonical transformation in Eq.~(\ref{tdir}) to the perturbation term in Eq.~(\ref{H1})  yields
\begin{equation}
\mathcal{H}_1\circ\mathcal{T}\equiv\mathcal{K}_1=-\frac{\mu}{B}\frac{1}{\rho},
\end{equation}
where the non-dimensional function $\rho=\rho(\phi,q,\Phi,Q)$ is
\begin{equation} \label{ro}
\rho=\sqrt{ 1+3\cos^2\phi+4\sigma\sin\phi+8\chi\cos\phi+4\sigma^2+4\chi^2 },
\end{equation}
in which the abscissa $Q/\omega$ and ordinate $q$ of the center of the reference ellipse have been scaled by its semi-minor and semi-major axes, respectively. That is
\begin{eqnarray} \label{chi}
\chi &=& \chi(q,\Phi)\equiv{q}/A, \\ \label{sigma}
\sigma &=& \sigma(Q,\Phi)\equiv(Q/\omega)/B,
\end{eqnarray}
where $A$ and $B$ are functions of $\Phi$ given by Eqs.~(\ref{semimajor}) and (\ref{semiminor}), respectively. With these abbreviations, the quadratic Hamiltonian in Eq.~(\ref{quadratic}) is rewritten as
\begin{equation}
\mathcal{K}_0=\omega\Phi\left(1-3\sigma^2\right).
\end{equation}
\par

Then, the coordinate $q$ is no longer cyclic for perturbed motion, and, in consequence, its conjugate momentum $Q$ will not remain constant anymore ---yet its variation will be small. Therefore, the abscissa of the center of the ellipse $x_\mathrm{C}=2Q/\omega$ will undergo small displacements.
\par

In the particular case of periodic unperturbed motion, the unperturbed value of $Q$ is zero, and hence $\sigma=0$. Then, $|\sigma|\ll1$ in the perturbed case. On the contrary, while the rate of variation of the ordinate of the center of the reference ellipse will suffer small variations due to the nonlinearities of the perturbed problem, the value $y_\mathrm{C}=q$ could be almost as large as the semi-major axis of the reference ellipse, and hence $\chi$ does not need to be small in itself.
\par

Integration of the perturbed problem $\mathcal{K}=\mathcal{K}_0+\mathcal{K}_1$ can be achieved in the averaging assumption that $\phi$ evolves much faster than the other variables, which evolve slowly, yet this solution have been claimed of doubtful practical application because it depends on special functions in a rather involved way \cite{Benest1976}. However, in view of the smallness of $\sigma$ one always can replace Eq.~(\ref{ro}) by the expansion
\begin{equation} \label{benest}
\rho=\sqrt{ 1+3\cos^2\phi+8\chi\cos\phi+4\chi^2 }+\mathcal{O}(\sigma).
\end{equation}
Then, a simpler approximation of the solution of the perturbed problem is obtained after neglecting terms of the order of $\sigma$ from the perturbation term (see \cite{Benest1976} for further details).
\par

In spite of the radical simplification obtained when neglecting $\mathcal{O}(\sigma)$ from Eq.~(\ref{benest}), the solution of the perturbation problem still depends on elliptic integrals, a fact that complicates improving Benest's solution by including higher order effects of $\sigma$, cf.~\cite{LidovVashkovyak1993}. Note, however, that, in the case of libration orbits about the primary, $\chi$ remains bounded and lower than the unity. For these orbits, Eq.~(\ref{benest}) can be further simplified by making the assumption $\chi=\mathcal{O}(\sigma^{m/n})$, with $n\ge{m}$. Choosing $m=n$ may be a realistic assumption for orbits with very small librations, but would be inaccurate otherwise. The assumption $\sigma=\mathcal{O}(\chi^2)$ is made here trying to broaden the applicability of the perturbation solution to encompass the case of orbits with larger librations.
\par

Under these assumptions, expansion of the inverse of $\rho$ in power series of a small parameter proportional to $\chi$ yields
\[
\frac{1}{\rho}=\sum_{n\ge0}\frac{1}{{\Delta}^{2n+1}}S_n,
\]
where
\begin{equation} \label{Delta}
\Delta=\sqrt{1+3\cos^2\phi},
\end{equation}
and the first $S_n$ coefficients are
\begin{eqnarray*}
S_0 &=& 1, \\
S_1 &=& -4\chi\cos\phi, \\
S_2 &=& -\frac{7}{2}\sigma\sin\phi-\frac{3}{2}\sigma\sin3\phi+7\chi^2+9\chi^2\cos2\phi, \\
S_3 &=& 30\sigma\chi\sin2\phi+9\sigma\chi\sin4\phi \\
&& -42\chi^3\cos\phi-22\chi^3\cos3\phi, \\[1ex]
S_4 &=& \frac{297}{4}\chi^4-\frac{149}{4}\sigma^2+\left(125\chi^4-\frac{501}{8}\sigma^2\right)\cos2\phi \\
&& +\left(\frac{227}{4}\chi^4-\frac{99}{4}\sigma^2\right)\cos4\phi -\frac{27}{8}\sigma^2\cos6\phi \\ 
&& -\frac{213}{2}\sigma\chi^2\sin\phi-\frac{627}{4}\sigma\chi^2\sin3\phi-\frac{153}{4}\sigma\chi^2\sin5\phi, \\
S_5 &=& \left(741\sigma^2\chi-495\chi^5\right)\cos\phi \\
&& +\left(561\sigma^2\chi-\frac{755}{2}\chi^5\right)\cos3\phi \\
&& +\left(207\sigma^2\chi-\frac{303}{2}\chi^5\right)\cos5\phi +27\sigma^2\chi\cos7\phi\\
&& +\frac{1585}{2}\sigma\chi^3\sin2\phi+670\sigma\chi^3\sin4\phi+\frac{285}{2}\sigma\chi^3\sin6\phi.
\end{eqnarray*}
\par

Besides, trying to better account for the effect of the central body's gravitation, the Hill problem Hamiltonian is reorganized in the new variables in the form of a perturbed harmonic oscillator in which both disturbing effects, viz.~the linear displacement of the center of the ellipse from the $x$ axis, on the one hand, and the gravitational attraction of the primary, on the other, are taken as perturbations of the same order.
\par

In the construction of the perturbation solution we found convenient to choose the following Hamiltonian arrangement: 
\begin{equation} \label{originalHam}
\mathcal{K}=\sum_{m\ge0}\frac{\epsilon^m}{m!}K_{m,0}(\phi,q,\Phi,Q),
\end{equation}
where $\epsilon$ is a formal small parameter used to indicate the importance of each perturbation term,
\begin{eqnarray} \label{lowoham}
K_{0,0} &=& \omega\Phi, \\ 
K_{m,0} &=& 0, \qquad 0<m<4, \\
K_{4,0} &=& -4!\left(3\omega^2B^2\sigma^2+\frac{\mu}{B}\frac{S_1}{\Delta}\right), \\ \label{highoham}
K_{m,0} &=& -\frac{\mu}{B}\frac{m!}{\Delta^{2m^*+1}}S_{m^*}, \quad m^*=m-4, \; m>4.
\end{eqnarray}
\par

The perturbation Hamiltonian given by Eqs.~(\ref{originalHam})--(\ref{highoham}), which is of two degrees of freedom, is reduced to a one degree of freedom Hamiltonian by the method of Lie-transforms \cite{Hori1966,Deprit1969}. Full details on the procedure, which is customary these days, can be consulted in standard textbooks like, for instance, \cite{Nayfeh2004,MeyerHall1992,BoccalettiPucacco1998v2}.
\par

Thus, we use the Lie-transforms method to compute an almost identity, canonical transformation 
\begin{equation} \label{toprimes}
\mathcal{T}:(\phi,q,\Phi,Q)\longrightarrow(\phi',q',\Phi',Q';\epsilon),
\end{equation}
from the original variables to new, prime canonical variables, such that it transforms Eq.~(\ref{originalHam}) into
\begin{equation}\label{newHam}
\mathcal{K}'=\mathcal{K}\circ\mathcal{T}\equiv\sum_{m=0}^{M}\frac{\epsilon^m}{m!}K_{0,m}(-,q',\Phi',Q')+\mathcal{O}(\epsilon^{M+1}).
\end{equation}
Then, after neglecting terms of the order of $\epsilon^{M+1}$ and higher, the coordinate $\phi'$ becomes cyclic, and, in consequence, $\Phi'$ is transformed into a formal integral. In this way, we obtain an integrable Hamiltonian of one degree of freedom in the new variables $q'$ and $Q'$.
\par

\section{Distant retrograde orbits: low-order solution} \label{s:low}

In the case of typical, almost-periodic, distant retrograde orbits, the massless body evolves far away from the origin at any time. Then, the gravitational attraction of the smaller primary is always a very small perturbation and, in consequence, a low-order perturbation solution captures the dynamics of the perturbed motion correctly.

Thus, neglecting terms of the order of $\epsilon^7$ and higher in Eq.~(\ref{newHam}), the normalized Hamiltonian becomes
\begin{eqnarray} \label{Ham60}
\mathcal{K}' &=& \omega\Phi'\bigg(1-3\sigma^2-2\frac{K(3/4)}{\pi}\frac{\mu}{\omega^2B^3} \\ \nonumber
&& 
-\frac{K(3/4)-E({3}/{4})}{(3/4)\pi}\frac{\mu}{\omega^2B^3}\chi^2 \bigg),
\end{eqnarray}
where $K(k^2)$ and $E(k^2)$ denote the complete elliptic integrals of the first and second kind, respectively, of modulus $k$. Note that the state functions $B$, $\sigma$, and $\chi$ are now functions of the new, prime variables, as given by Eqs.~(\ref{semiminor}), (\ref{chi}) and (\ref{sigma}), respectively. Namely, $B=\sqrt{2\Phi'/\omega}$, $\chi=\frac{1}{2}q'/B$, and $\sigma=(Q'/\omega)/B$. Besides, to shorten notation, the abbreviations
\[
\tilde{K}\equiv\frac{1}{\pi}K({3}/{4}), \qquad
\tilde{E}\equiv\frac{1}{\pi}E({3}/{4}),
\]
are used in what follows. Finally, we find convenient to define the libration frequency
\begin{equation} \label{Om}
\Omega=\Omega(\Phi')\equiv\sqrt{(\tilde{K}-\tilde{E})\frac{\mu}{B^3}}.
\end{equation}
\par

Then, Eq.~(\ref{Ham60}) is rewritten
\begin{eqnarray} \label{Ham6}
\mathcal{K}' &=& \omega\Phi'\bigg[1-3\sigma^2
-\frac{\Omega^2}{\omega^2}\left(\frac{2\tilde{K}}{\tilde{K}-\tilde{E}}+\frac{4}{3}\chi^2\right) \bigg],
\end{eqnarray}
in which the new momentum $\Phi'$ is constant because its conjugate variable $\phi'$ has been removed up to the truncation order. This formal integral decouples the Hamiltonian flow derived from Eq.~(\ref{Ham6}) into the reduced system
\begin{eqnarray} \label{qp7}
\frac{\mathrm{d}q'}{\mathrm{d}t} &=& -3Q', \\  \label{Qp7}
\frac{\mathrm{d}Q'}{\mathrm{d}t} &=& \frac{1}{3}\Omega^2q',
\end{eqnarray}
and the differential equation
\begin{equation} \label{fp7}
\frac{\mathrm{d}\phi'}{\mathrm{d}t}=\frac{\partial\mathcal{K}'(q'(t),Q(t);\Phi')}{\partial\Phi'},
\end{equation}
which is integrated by quadrature after solving Eqs.~(\ref{qp7})--(\ref{Qp7}).
\par

The solution of Eqs.~(\ref{qp7})--(\ref{fp7}) is standard, yielding
\begin{eqnarray} \label{solq6}
q' &=& q'_*\cos\Omega{t}-p_*\sin\Omega{t}, \\ \label{solQ6}
Q' &=& Q'_*\cos\Omega{t}+(q'_*\Omega/3)\sin\Omega{t}, \\ \label{solf6}
\phi' &=& \phi'_*+\omega\left(1+\delta\right)t
+\frac{\Omega}{\omega}\frac{{q'_*}^2-p_*^2}{8B^2}\sin2\Omega{t} \\ \nonumber
&& 
+\frac{\Omega}{\omega}\frac{q'_*p_*}{4B^2}(\cos2\Omega{t}-1),
\end{eqnarray}
where we abbreviated $p_*=3Q'_*/\Omega$ and
\begin{equation} \label{delta7}
\delta=\left(\frac{\tilde{K}}{\tilde{K}-\tilde{E}}+\frac{{q'_*}^2+p_*^2}{4B^2}\right)\frac{\Omega^2}{\omega^2}.
\end{equation}
That is, on average, the center of the reference ellipse evolves in the $Q,q$ plane with harmonic oscillations of frequency $\Omega$. On the other hand, the linear growing of the phase $\phi'$ with time, which happens at the perturbed rate $\omega(1+\delta)$, is modulated with periodic oscillations of frequency $2\Omega$.
\par

Orbit and libration periods, which are noted $T$ and $T^*$ respectively, are defined as
\begin{eqnarray} \label{periodo6}
T &=& \frac{2\pi}{\omega(1+\delta)}, \\ \label{periodl6}
T^* &=& \frac{2\pi}{\Omega}.
\end{eqnarray}
Note that Eq.~(\ref{periodl6}), with $\Omega$ given by Eq.~(\ref{Om}), matches exactly the expression of the limiting period given by Benest for the case of librations of small amplitude, cf.~Eq.~(22) of \cite{Benest1976}. Remarkably, the libration period still exists for periodic orbits centered at the origin, for which $q'_*=0$, $Q'_*=0$. In that case, $\phi'$ grows linearly with perturbed frequency $\omega[1+\tilde{K}\mu/(\omega^2B^3)]$.
\par

The solution given by Eqs.~(\ref{solq6})--(\ref{solf6}), in prime variables, must be complemented with the periodic corrections in Appendix \ref{a:LT}, up to the sixth order, that define the transformation from prime to original variables. However, since these corrections are generally small, we will see that the prime variables solution, which dominates the long-term dynamics, may provide by itself a good description of the motion.
\par

Sample illustrations for different kinds of motion are presented below, where we use Hill problem units, in which $\mu=1$ and $\omega=1$. Results provided by the analytical approximation are compared with the true orbit as given by the numerical integration of the original Hamiltonian flow in Eqs.~(\ref{H0})--(\ref{H1}). Namely,
\begin{eqnarray} \label{xp}
\frac{\mathrm{d}x}{\mathrm{d}t} &=& X+y, \\
\frac{\mathrm{d}y}{\mathrm{d}t} &=& Y-x, \\
\frac{\mathrm{d}X}{\mathrm{d}t} &=& -\frac{x}{r^3}+2x+Y, \\ \label{Yp}
\frac{\mathrm{d}Y}{\mathrm{d}t} &=& -\frac{y}{r^3}-y-X.
\end{eqnarray}

Three different test cases are discussed. The first one corresponds to a distant retrograde periodic orbit that would correspond to a true periodic orbit of the original Hill problem. The second case deals with an almost periodic orbit with small libration amplitude. Finally, the third case examines a quasi-satellite orbit with large librations about the smaller primary. Corresponding initial conditions of the test orbits are provided in Table \ref{t:test6}, which also includes the theoretical values of the orbital and librational periods predicted by Eqs.~(\ref{periodo6}) and (\ref{periodl6}), respectively.
\par

\begin{table*}[htbp]
\centering
\begin{tabular}{@{}lrrrllll@{}} 
\multicolumn{1}{c}{Case} & \multicolumn{1}{c}{$x_*$} & \multicolumn{1}{c}{$y_*$} & \multicolumn{1}{c}{$X_*$} & \multicolumn{1}{c}{$Y_*$} & \multicolumn{1}{c}{$T$} & \multicolumn{1}{c}{$T^*$} & \multicolumn{1}{c}{$\Phi$}  \\
\noalign{\smallskip}\hline\noalign{\smallskip}
1: periodic        & $0.1$ & $20.0$ & $-10.0$ & $-0.1$ & $6.27888$ & $362.215$ & $50.005$ \\
2: small amplitude libration  & $0.1$ & $20.0$ & $-10.5$ & $-0.1$ & $6.27815$ & $335.394$ & $45.130$ \\
3: large amplitude libration & $0.0$ & $10.0$ &  $-0.5$ & $-0.1$ & $6.27611$ & $335.477$ & $45.145$ \\
\noalign{\smallskip}\hline
\end{tabular}
\caption{Test cases for the sixth order solution}
\label{t:test6}
\end{table*}

\subsection{Case 1}
In case 1, the initial conditions chosen are $x_*=0.1$, $y_*=20$, $X_*=-10$, $Y_*=-0.1$, fulfilling the conditions for a periodic orbit centered at the origin in the unperturbed problem approximation given by Eqs.~(\ref{periodic0}) and (\ref{centered0}). Then, the values $q_*=Q_*=0$, $\phi_*=0.00999967$, and $\Phi=50.005$, are obtained from Eq.~(\ref{tinv}). The analytical solution is formulated in the prime variables, so the inverse transformation in Appendix \ref{a:LT} should be applied first, up to the order 6, to compute the correct initial conditions to be used when evaluating Eqs.~(\ref{solq6})--(\ref{solf6}). However, due to the low-order of the analytical solution, the effects of using these corrections are not too much relevant in what respects to the average behavior, and the values of $q_*$, $Q_*$, $\phi_*$, and $\Phi$, given above are used directly as if they were prime elements. The orbit is then propagated using the solution in Eqs.~(\ref{solq6})--(\ref{solf6}) for a libration period $T^*=362.215$, as given by Eq.~(\ref{periodl6}), and the analytical results are then compared with the true solution for the same initial conditions, which is obtained from the numerical integration of Eqs.~(\ref{xp})--(\ref{Yp}).
\par

\begin{figure}[htbp]
\begin{center}
\includegraphics[scale=0.8, angle=0]{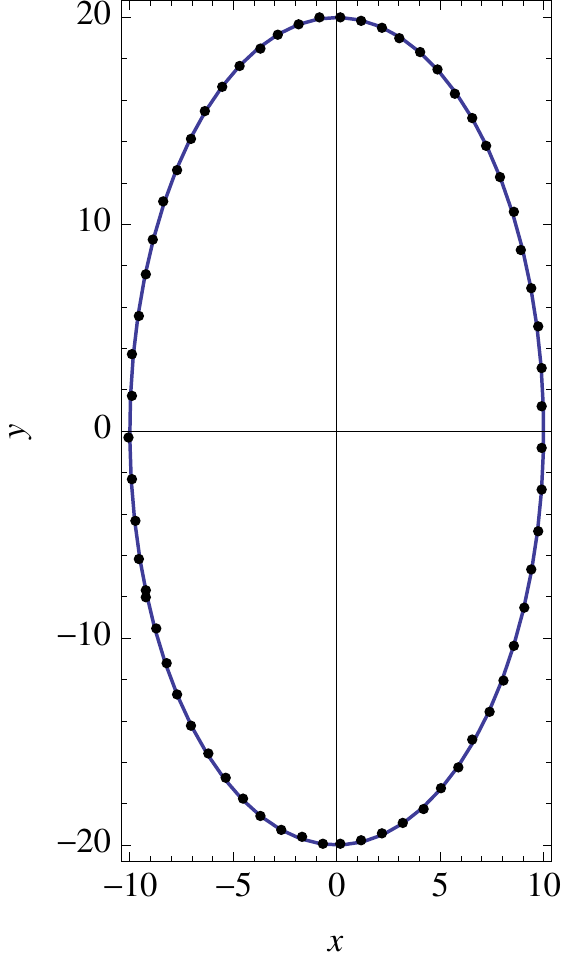}
\caption{Case 1 of Table \protect\ref{t:test6}. Last orbital periods of the true orbit (black dots) and the low-order analytical solution (full line).}
\label{f:orbit1six}
\end{center}
\end{figure}

\begin{figure}[htb]
\begin{center}
\includegraphics[scale=0.7, angle=0]{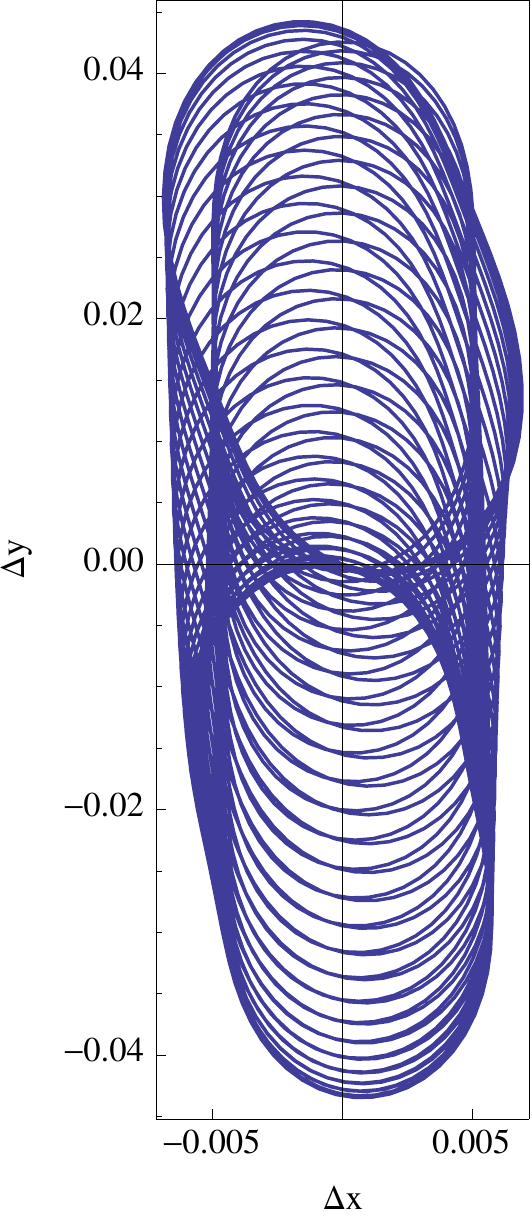}\quad
\includegraphics[scale=0.7, angle=0]{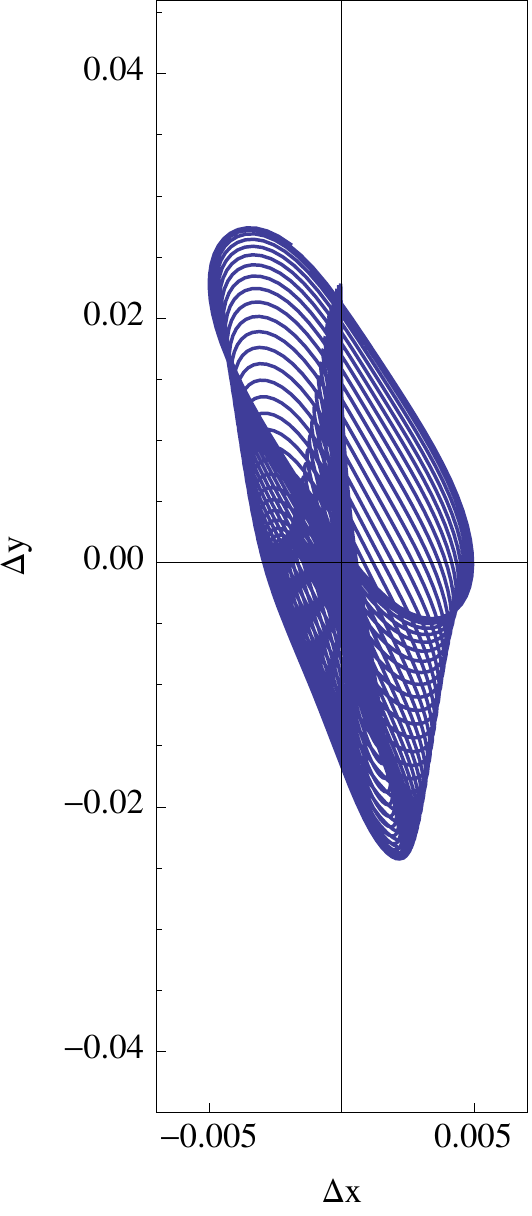}
\caption{Case 1 of Table \protect\ref{t:test6}. Errors of the analytical solution. Left: averaged equations only. Right: analytical approximation with short period corrections. Note the different scales of the abscissas and ordinates axes.}
\label{f:error1six}
\end{center}
\end{figure}

As shown in Fig.~\ref{f:orbit1six}, the analytical solution provides a good approximation of the true orbit, and, along one libration period, both orbits match at the precision of the graphics. The smallness of the errors between the analytical solution for the averaged motion given by the Eqs.~(\ref{solq6})--(\ref{solf6}) and the true orbit along one libration period, which are of the order of just a few thousandths relative to the size of the orbit, can be appreciated in the left plot of Fig.~\ref{f:error1six}. Using the inverse transformation to compute the correct initial conditions in prime variables, and recovering the short-period effects removed by the averaging at each selected point in which the analytical theory is evaluated, improves, of course, the accuracy of the analytical solution, as illustrated with the right plot of Fig.~\ref{f:error1six}. Nevertheless, these improvements are just moderate because of the early truncation of the perturbation theory. Analogous errors are obtained for the conjugate momenta, and are not displayed.
\par

In spite of the initial conditions correspond to an orbit that is periodic only on average, the true orbit is almost periodic, and, after the theoretical orbital period $T=6.27888$ predicted by Eq.~(\ref{periodo6}), we obtain a periodicity error, defined as the maximum of the relative errors of the position and velocity vectors
\begin{eqnarray*}
\varepsilon &=&
\max\left(\sqrt{\frac{(x(T)-x_*)^2+(y(T)-y_*)^2}{x_*^2+y_*^2}},\right. \\
&& \left.\sqrt{\frac{(X(T)-X_*)^2+(Y(T)-Y_*)^2}{X_*^2+Y_*^2}}\right),
\end{eqnarray*}
of the order of $10^{-5}$. 
\par

If desired, the case 1 initial conditions and period can be easily improved by differential corrections to get a true periodic orbit of the Hill problem. In particular, when using the algorithm in \cite{LaraPelaez2002} and working in double precision arithmetic, three iterative corrections are enough to converge to a ``true'' periodic orbit, with periodicity error $\varepsilon=\mathcal{O}(10^{-15})$. 
\par

\begin{figure}[htbp]
\begin{center}
\includegraphics[scale=0.7, angle=0]{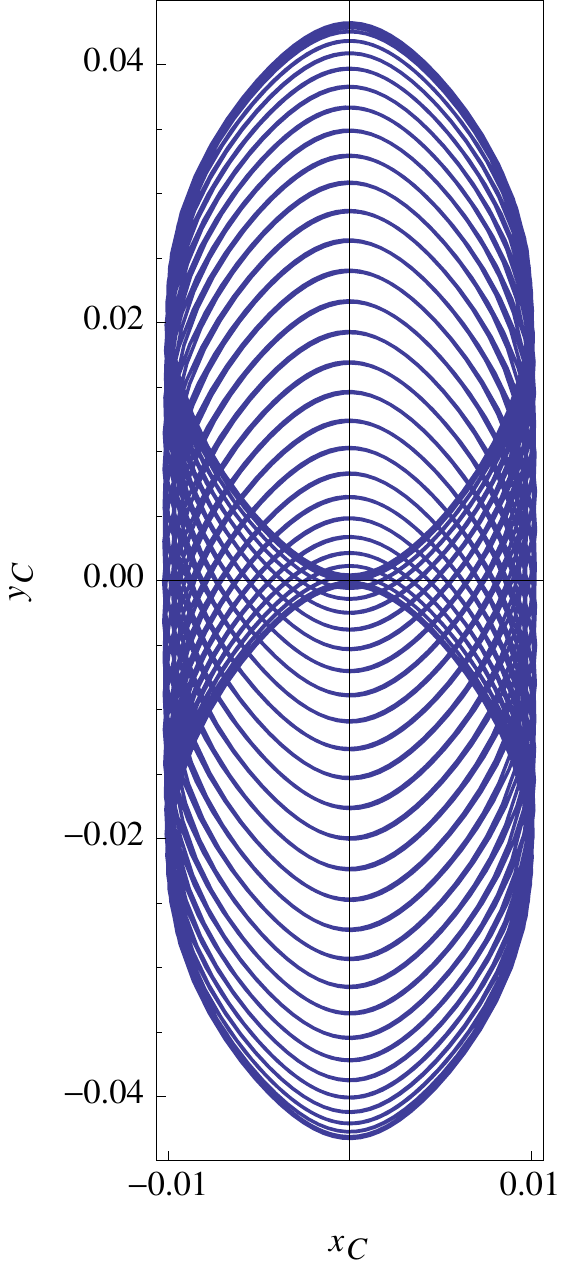} 
\includegraphics[scale=0.7, angle=0]{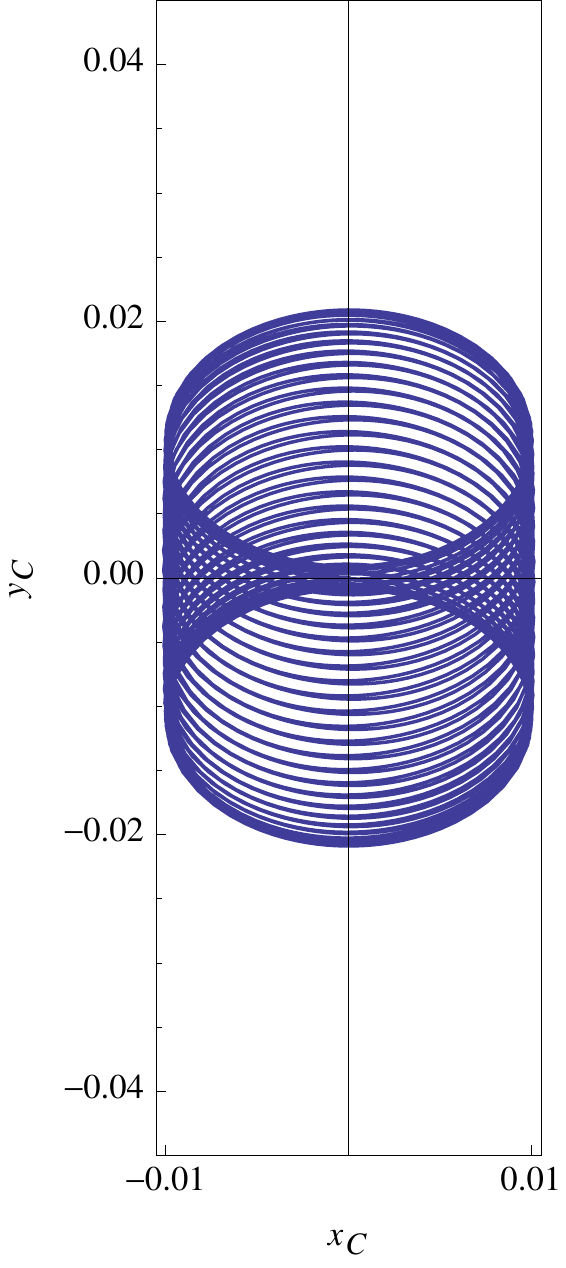} 
\caption{Case 1 of Table \protect\ref{t:test6}. Coordinates ($x_C,y_C$) of the center of the reference ellipse of the true solution (left plot) and analytical solution with periodic corrections (right plot).}
\label{f:orbit1sixce}
\end{center}
\end{figure}

On the other hand, the center of the reference ellipse only remains at the origin on average, but it is affected of small periodic effects which become apparent when the solution of Eqs.~(\ref{solq6})--(\ref{solf6}) is complemented with the short-period corrections. This is illustrated in Fig.~\ref{f:orbit1sixce}, where the left plot shows the true oscillations of the center of the reference ellipse as obtained from the numerical integration, and the right plot shows corresponding analytical predictions. By comparison of Figs.~\ref{f:orbit1six} and \ref{f:error1six} it becomes evident that the errors in the propagation of the orbit are mostly a consequence of the errors of the analytical theory in modeling the short-period motion of the center of the reference ellipse.
\par

\subsection{Case 2}

\begin{figure}[htb]
\begin{center}
\includegraphics[scale=0.8, angle=0]{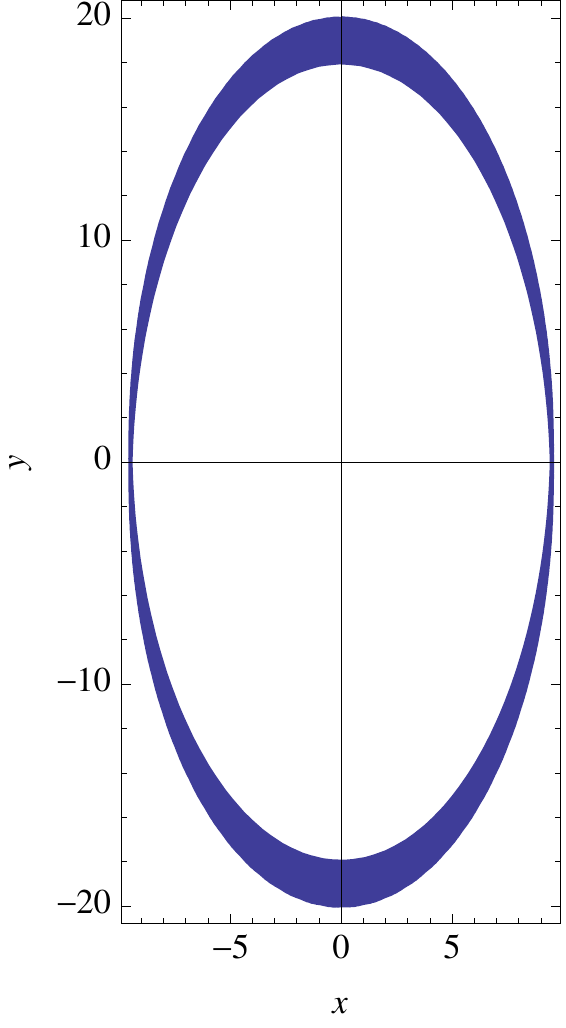}
\caption{Orbit with small libration corresponding to case 2 of Table \protect\ref{t:test6}.}
\label{f:orbit2six}
\end{center}
\end{figure}

When $y_*+2X_*\ne0$ the condition in Eq.~(\ref{centered0}) for a centered ellipse of the non-perturbed problem is broken. Therefore, $q_*\ne0$ and the reduced solution given by Eqs.~(\ref{solq6}) and (\ref{solQ6}) no longer vanishes. On the contrary, the center of the osculating ellipse of both kind of solutions, analytical and numerical,  experiences a librational motion of non-negligible amplitude along the $y$ axis. For the test case 2 of Table \ref{t:test6}, the effects of this libration are illustrated in Fig.~\ref{f:orbit2six}, where the analytical solution given by Eqs.~(\ref{solq6})--(\ref{solf6}) is superimposed to the numerical solution of Eqs.~(\ref{xp})--(\ref{Yp}), showing that both orbits match at the precision of the graphics. 
\par

\begin{figure}[htb]
\begin{center}
\includegraphics[scale=0.7, angle=0]{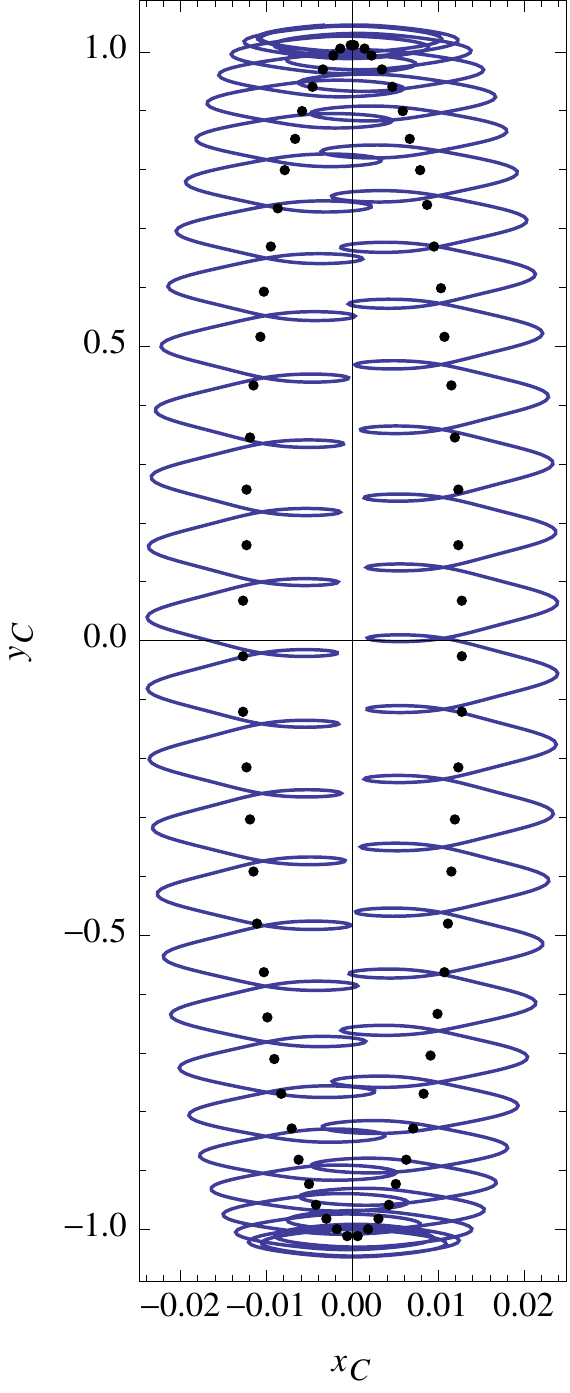}
\includegraphics[scale=0.7, angle=0]{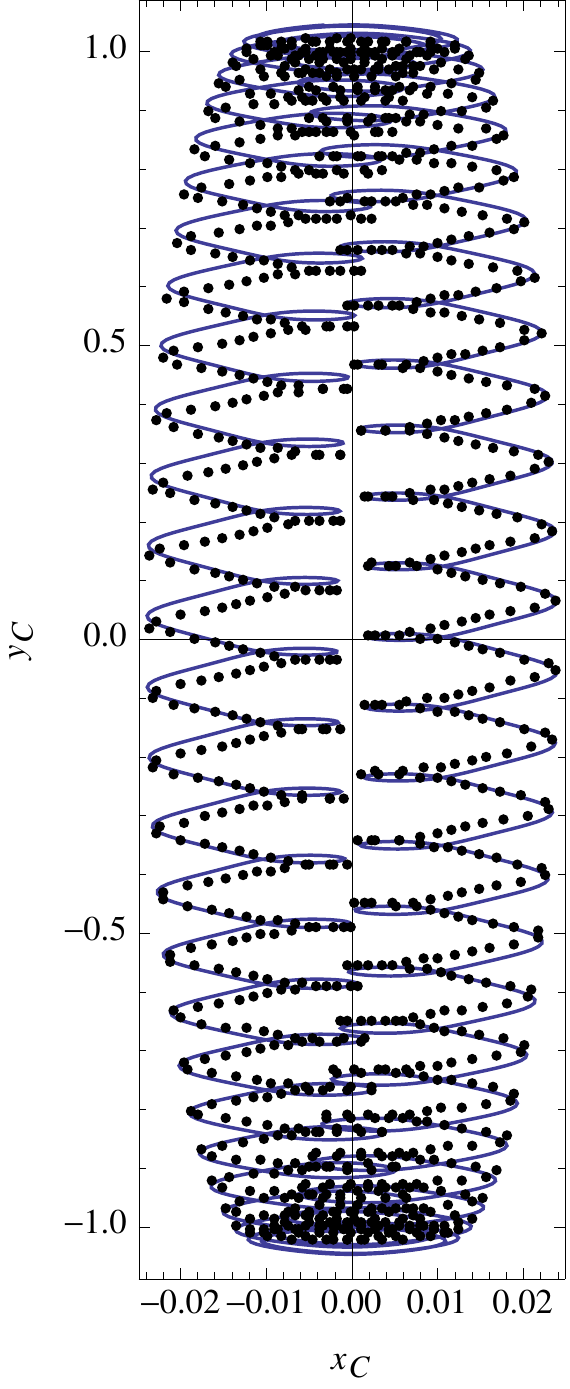}
\caption{Case 2 of Table \protect\ref{t:test6}. Trajectory $(x_C,y_C)$ of the center of the reference ellipse of the analytical solution (black dots) superimposed to the true solution (full line) along one libration period. Left: averaged terms only. Right: analytical solution with short-period corrections. Note the different scales of the $x_C$ and $y_C$ axes.}
\label{f:orbit2sixce}
\end{center}
\end{figure}

The libration effects are clearly observed in the motion of the center of the reference ellipse, whose amplitude is now clearly larger than in the previous example, as observed when comparing Fig.~\ref{f:orbit2sixce} with previous Fig.~\ref{f:orbit1sixce}. Again, the averaged solution obtained by evaluation of Eqs.~(\ref{solq6}) and (\ref{solQ6}) only differs in periodic effects from its numerically integrated equivalent, as shown in the left plot of Fig.~\ref{f:orbit2sixce}, thus correctly predicting the average motion of the actual libration. The orbit is better approximated when the short-period effects removed by the averaging are incorporated to the solution, a procedure that includes applying the inverse transformation to compute the correct initial conditions in prime variables. Still, as shown in the right plot of Fig.~\ref{f:orbit2sixce}, there is a shift between the analytical approximation and the numerical solution due to the low-order truncation of the perturbation approach.
\par

The errors between the analytical solution and the true orbit of the test case 2 remain small, of similar magnitude as those showed in Fig.~\ref{f:error1six} for test case 1, and hence are not displayed. 
\par

\subsection{Case 3}

The low-order analytical solution still remains valid for  distant retrograde orbits with large librations, yet, as expected, the accuracy of the predictions degrade. Indeed, the large amplitude of the librational motion makes that, eventually, the body of negligible mass gets much closer to the smaller primary than in previous examples. During these close passages, the strength of the gravitational effect of the primary is undervalued by the analytical solution due to the early truncation of the perturbation approach. Because of that, the low-order analytical solution predicts larger librations than actually happen. However, even though the quantitative differences between the true and predicted orbits may be significant, the qualitative behavior still remains quite similar. 
\par

\begin{figure}[htb]
\begin{center}
\includegraphics[scale=0.8, angle=0]{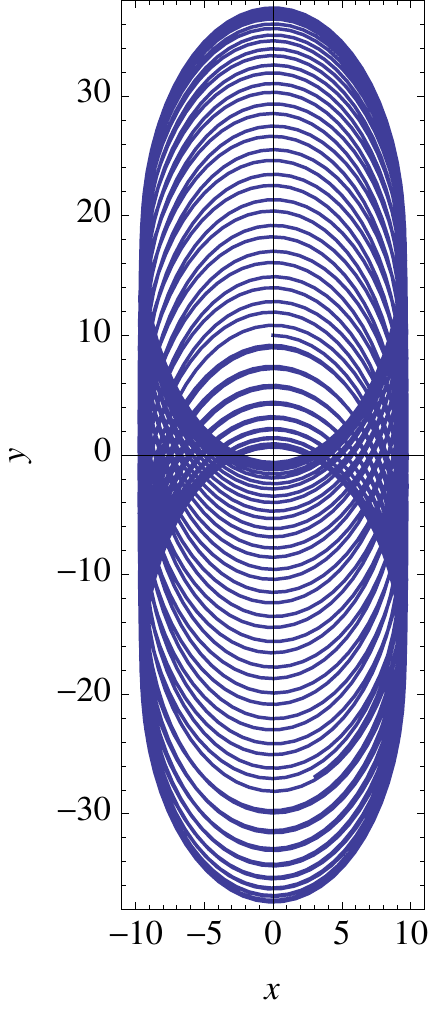} \qquad
\includegraphics[scale=0.8, angle=0]{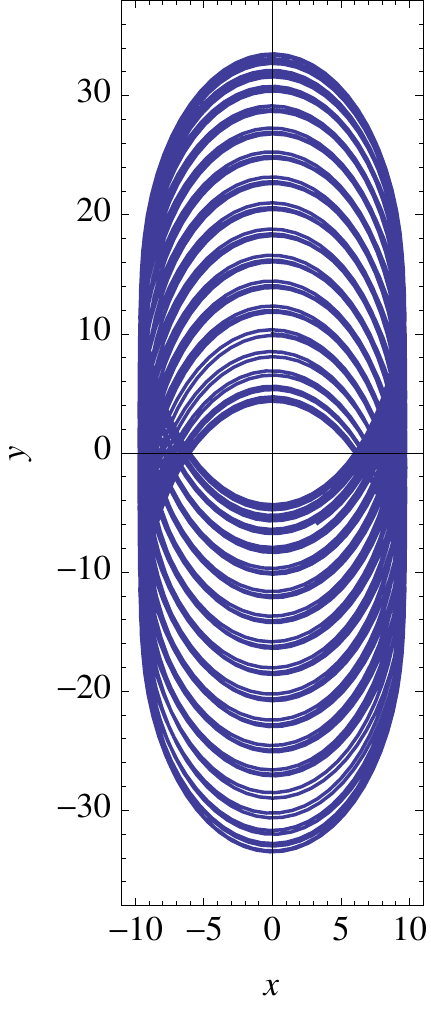}
\caption{Case 3 of Table \protect\ref{t:test6}. Left: analytical orbit. Right: actual orbit.}
\label{f:orbit3six}
\end{center}
\end{figure}

\begin{figure}[htb]
\begin{center}
\includegraphics[scale=0.7, angle=0]{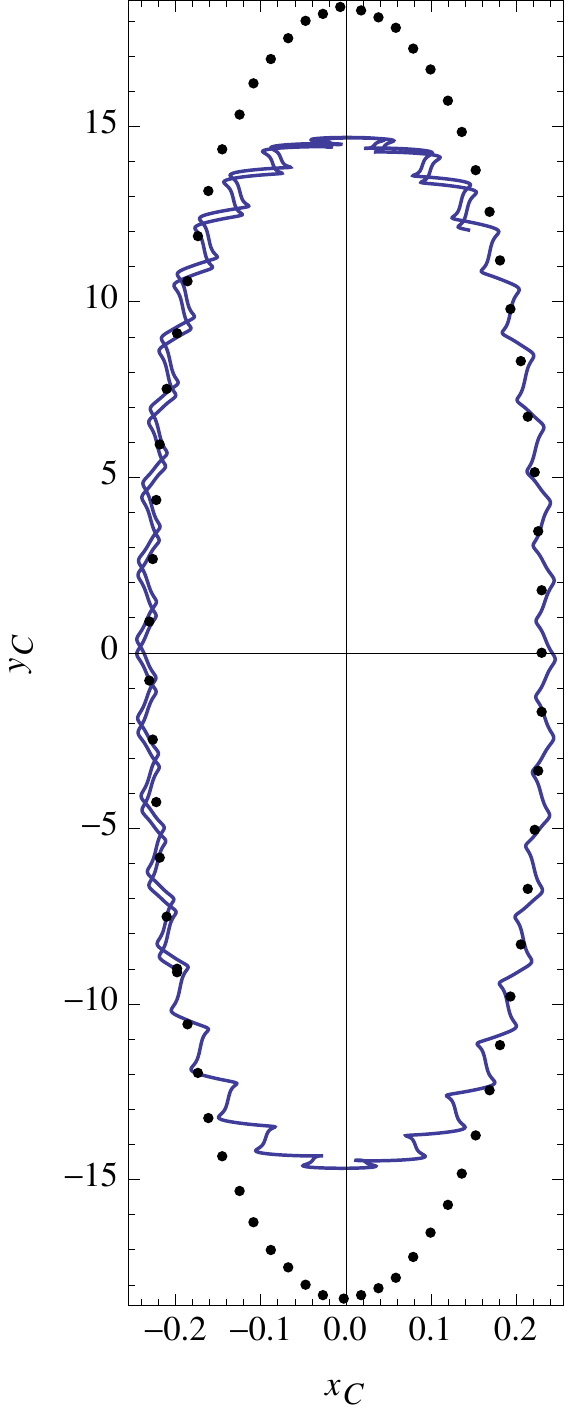}
\caption{Case 3 of Table \protect\ref{t:test6}. Trajectory $(x_C,y_C)$ of the center of the reference ellipse of the analytical solution (black dots) superimposed to the true solution (full line) along one libration period. Note the different scales of the $x_C$ and $y_C$ axes.}
\label{f:orbit3sixce}
\end{center}
\end{figure}

The effects of the mismodeling of the dynamics in orbits with large librations are illustrated in Figs.~\ref{f:orbit3six} and \ref{f:orbit3sixce} for the test case 3 of Table \ref{t:test6}. Thus, while the analytical solution predicts an amplitude of the librational motion along the $y$ axis direction of about $\pm18$ length units, which is travelled in a libration period $T^*=335.48$, the amplitude of the true libration is notably smaller, of only about $\pm14$ length units, which are travelled in the much shorter libration period of about $232$ time units. In consequence, the orbit predicted by the analytical solution approaches much closer to the primary than the real orbit. 
Indeed, as shown in the left plot of Fig.~\ref{f:orbit3six}, the incorrect modeling of the smaller primary gravitation due to the low-order truncation of the perturbation approach, makes that the amplitude of the librations predicted by the analytical solution would allow the body of negligible mass to enter the Hill region. Quite on the contrary, the actual dynamics prevents the body of negligible mass to get closer to the origin than about 6 times the Hill radius, as shown in the right plot of Fig.~\ref{f:orbit3six}.
\par

For this particular example, we note in Fig.~\ref{f:orbit3sixce} that the actual motion of the center of the ellipse is almost periodic after a libration period ---that we checked that is of $232.04$ time units. If desired, the periodicity of the test case 3 orbit is easily improved by differential corrections leading to a slightly unstable, true periodic orbit ($\varepsilon=\mathcal{O}(10^{-13})$) with initial conditions
\[
\begin{array}{rcl}
x &=& \phantom{-}0.0009558942643146, \\
y &=& \;10.09070684586246, \\
X &=&-0.5908147794362844, \\
Y &=&-0.1003142256682326,
\end{array}
\]
and period $T=232.2079125513217$.
\par

\section{Quasi-satellite orbits with large libration} \label{s:high}

The applicability of the analytical approach can be extended by computing higher orders of the perturbation theory. However, in spite of the integrability of the higher order averaged Hamiltonian, the benefits of having an exact, explicit solution as the one in Eqs.~(\ref{solq6})--(\ref{solf6}) are immediately lost.
\par

Indeed, the next non-vanishing term of the new Hamiltonian, which is the 8th order term in our Hamiltonian arrangement, includes the fourth power of $q'$. Thus,
\begin{equation} \label{Ham8}
\mathcal{K}'=\omega\Phi'(b_1-b_2\sigma^2-2b_3\chi^2-b_4\chi^4),
\end{equation}
where the non-dimensional functions $b_i=b_i(\Phi')$, $i=1,\dots4$, are
\begin{eqnarray*}
b_1 &=& 1-\frac{2}{\tilde{K}-\tilde{E}}\left(\tilde{K}
+\frac{\tilde{K}^2-1/2}{\tilde{K}-\tilde{E}}\frac{\Omega^2}{\omega^2}\right)\frac{\Omega^2}{\omega^2}, \\
b_2 &=& 3\left(1+\frac{4}{9}\frac{4\tilde{E}-\tilde{K}}{\tilde{K}-\tilde{E}}\frac{\Omega^2}{\omega^2}\right), \\
b_3 &=& \frac{2}{3}\frac{\Omega^2}{\omega^2}, \\
b_4 &=& \frac{1}{9}\frac{11\tilde{K}-14\tilde{E}}{\tilde{K}-\tilde{E}}\frac{\Omega^2}{\omega^2},
\end{eqnarray*}
where $b_2$, $b_3$ and $b_4$ are strictly positive quantities. On the contrary, the sign of $b_1$ may change depending on the value of $\Omega^2/\omega^2$. Indeed, the condition $b_1=0$ leads to
\[
\frac{\Omega^2}{\omega^2}=(\tilde{K}-\tilde{E})\frac{\tilde{K}\pm\sqrt{3\tilde{K}^2-1}}{2\tilde{K}^2-1}.
\]
Besides, $0<\Omega^2\ll\omega^2$ for a perturbation problem. Hence, $b_1>0$ when
\[
\frac{\Omega^2}{\omega^2}\le(\tilde{K}-\tilde{E})\frac{\tilde{K}-\sqrt{3\tilde{K}^2-1}}{2\tilde{K}^2-1}\approx0.226
\]
which will be the typical case of distant retrograde orbits.
\par

From brevity, in what follows we suppress the prime notation provided there is no risk of confusion. Then, from the Hamilton equations of Eq.~(\ref{Ham8}),
\begin{equation} \label{dq8}
\frac{\mathrm{d}q}{\mathrm{d}t}=\frac{\partial\mathcal{K}}{\partial{Q}}=-b_2Q.
\end{equation}

Besides, for any energy level $\mathcal{K}=h$, $\sigma$, and, therefore, $Q$, are trivially solved from Eq.~(\ref{Ham8}), and can be replaced in Eq.~(\ref{dq8}) to give a differential equation in separate variables that is solved by quadrature. However, this analytical solution does not help much, because it provides $t$ as a function of $q$, an implicit solution that must be inverted to compute $q$ and $Q$ as functions of time so that $\phi$ can be solved next from its Hamilton equation. In addition, the quadrature solving Eq.~(\ref{dq8}) involves the incomplete elliptic integral of the first kind, a fact that introduces supplementary complications in the process of expressing $q$, $Q$, and $\phi$ like explicit functions of time.
\par

\begin{figure}[htb]
\begin{center}
\includegraphics[scale=0.8, angle=0]{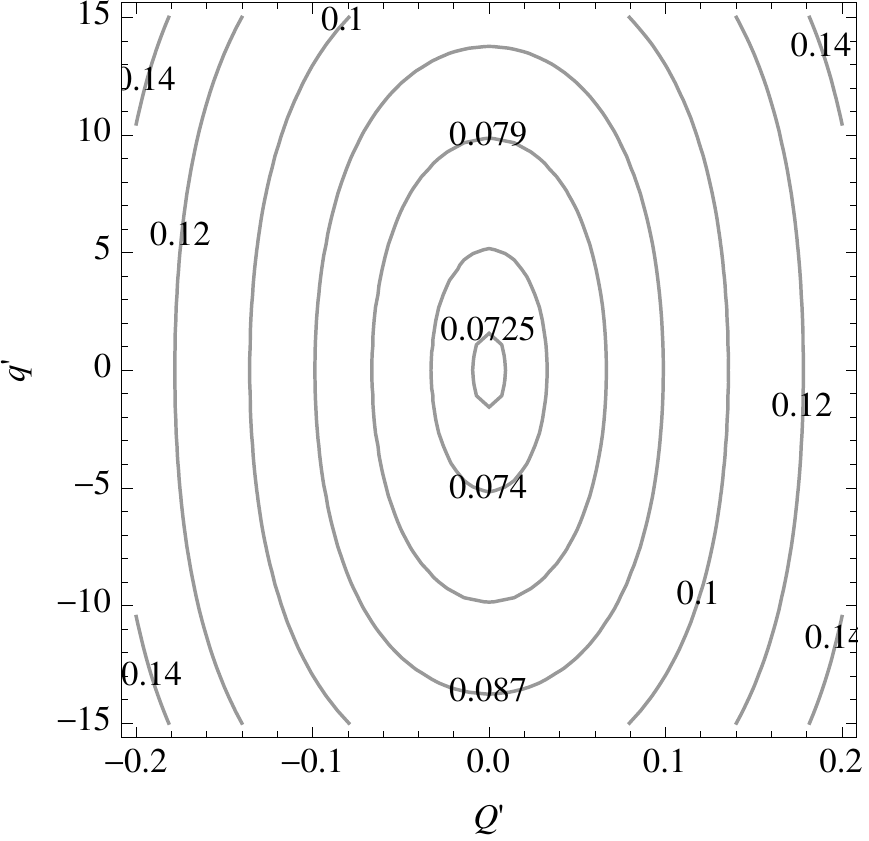}
\caption{Contours $\mathcal{K}'=h$ of the eight order averaged Hamiltonian in Eq.~(\protect\ref{Ham8}) for $\Phi'=45$ ($\mu=\omega=1$).}
\label{f:contours8F45}
\end{center}
\end{figure}

On the other hand,  simple contour plots of the Hamiltonian (\ref{Ham8}), which are computed without need of carrying out any integration, show that the orbits of the reduced phase space remain in this higher order approximation like periodic perturbed harmonic oscillations, as illustrated in Fig.~\ref{f:contours8F45}. Hence, an approximate solution to the differential system comprised by Eq.~(\ref{dq8}) and 
\begin{equation} \label{dQ8}
\frac{\mathrm{d}Q}{\mathrm{d}t}=-\frac{\partial\mathcal{K}}{\partial{q}}
=-\frac{\omega^2}{4}q\left(2b_3+\frac{\omega{b}_4}{4\Phi}q^2\right),
\end{equation}
can be computed using Lindstedt-Poincar\'e method.
\par

Thus, we improve the low-order solution in Eqs.~(\ref{solq6}) and (\ref{solQ6}) by making first the change of the independent variable
\begin{equation}
\tau=nt.
\end{equation}
Then, the differential system in Eqs.~(\ref{dq8}) and (\ref{dQ8}) is rewritten in the new independent variable as
\begin{eqnarray} \label{qp8tau}
n\frac{\mathrm{d}q(\tau)}{\mathrm{d}\tau} &=& -b_2Q(\tau), \\ \label{Qp8tau}
n\frac{\mathrm{d}Q(\tau)}{\mathrm{d}\tau} &=&
-\frac{1}{4}\omega^2q(\tau)\left[2b_3+\frac{{b}_4}{2}\frac{q(\tau)^2}{B^2}\right].
\end{eqnarray}
Next, $n$, $q$, and $Q$ are replaced by the series
\begin{eqnarray} \label{sn}
n &=& 1+\sum_{i\ge1}\epsilon^in_i, \\ \label{sq}
q &=& \sum_{i\ge0}\epsilon^iq_i(\tau), \\ \label{sQ}
Q &=& \sum_{i\ge0}\epsilon^iQ_i(\tau),
\end{eqnarray}
where, once more, $\epsilon$ is a formal small parameter used to indicate the relative strength of each term of the respective series. It follows the replacement of Eqs.~(\ref{sn})--(\ref{sQ}) into the new differential system given by Eqs.~(\ref{qp8tau})--(\ref{Qp8tau}), in which the coefficients of the same powers of $\epsilon$ are then identified. In this way we obtain a chain of differential systems that can be solved sequentially. Finally, the coefficients $n_i$ are chosen in such a way that the solution of each differential system is periodic.
\par

The zeroth order terms of the series give the system
\begin{eqnarray} \label{qp0}
\frac{\mathrm{d}q_0}{\mathrm{d}\tau} &=& -3Q_0, \\ \label{Qp0}
\frac{\mathrm{d}Q_0}{\mathrm{d}\tau} &=& \frac{1}{3}\Omega^2q_0,
\end{eqnarray}
whose solution matches the 6th order perturbation solution in Eqs.~(\ref{solq6})--(\ref{solQ6}), but that must be evaluated in the new time scale, viz.
\begin{eqnarray} \label{solq0}
q_0 &=& q_*\cos\Omega\tau-p_*\sin\Omega\tau, \\ \label{solQ0}
Q_0 &=& Q_*\cos\Omega\tau+(q_*\Omega/3)\sin\Omega\tau.
\end{eqnarray}
\par

The next differential system is
\begin{eqnarray} \label{qp1}
\frac{\mathrm{d}q_1}{\mathrm{d}\tau} &=& -3Q_1+3\left(e_1\frac{\Omega^2}{\omega^2}-n_1\right)Q_0(\tau), \\ \label{Qp1}
\frac{\mathrm{d}Q_1}{\mathrm{d}\tau} &=& \frac{1}{3}\Omega^2q_1+\frac{1}{3}\Omega^2q_0(\tau)\left[e_2\frac{q_0^2(\tau)}{B^2}-n_1\right],
\end{eqnarray}
in which $q_0(\tau)$ and $Q_0(\tau)$ must be replaced by the right sides of Eqs.~(\ref{solq0}) and (\ref{solQ0}), respectively, and we shortened notation making
\begin{eqnarray*}
e_1 &=& \frac{4}{9}\frac{4\tilde{E}-\tilde{K}}{\tilde{K}-\tilde{E}}\approx1.26, \\ 
e_2 &=& \frac{1}{24}\frac{11\tilde{K}-14\tilde{E}}{\tilde{K}-\tilde{E}}\approx0.298.
\end{eqnarray*}
After solving Eqs.~(\ref{qp1})--(\ref{Qp1}), we choose
\begin{equation} \label{n1}
n_1=\frac{1}{2}e_1\frac{\Omega^2}{\omega^2}+\frac{3}{8}e_2\frac{q_*^2+p_*^2}{B^2},
\end{equation}
to get the periodic solution
\begin{eqnarray} \label{solq1}
q_1 &=& e_2\frac{q_*^2-3p_*^2}{32B^2}q_*(\cos3\Omega\tau-\cos\Omega\tau) \\ \nonumber
&& +\left(e_2\frac{21q_*^2+9p_*^2}{32B^2}-\frac{e_1}{2}\frac{\Omega^2}{\omega^2}\right)p_*\sin\Omega\tau \\ \nonumber
&& -e_2\frac{3q_*^2-p_*^2}{32B^2}p_*\sin3\Omega\tau \\
\label{solQ1}
Q_1 &=& 3e_2\frac{3q_*^2-p_*^2}{32B^2}Q_*(\cos3\Omega\tau-\cos\Omega\tau) \\  \nonumber
&& +\frac{1}{6}\left(e_2\frac{11q_*^2+15p_*^2}{16 B^2}-e_1\frac{\Omega^2}{\omega^2}\right)q_*\Omega\sin\Omega\tau \\  \nonumber
&& +e_2\frac{q_*^2-3p_*^2}{32B^2}q_*\Omega\sin3\Omega\tau.
\end{eqnarray}
\par

The higher order solution 
\begin{eqnarray} \label{qLin}
q=q_0(\tau)+q_1(\tau), \\ \label{QLin}
Q=Q_0(\tau)+Q_1(\tau),
\end{eqnarray}
which must be evaluated in the time argument
\begin{equation} \label{tau}
\tau=(1+n_1)t,
\end{equation}
allows to solve $\phi$ from the corresponding Hamilton equation derived from the Eq.~(\ref{Ham8}). Besides, the short-period effects removed by the averaging are recovered using the transformation equations from prime to original variables up to the 8th order of the perturbation approach. In addition, a refined estimation of the libration period is given by
\begin{equation} \label{periodl8}
T^*=\frac{2\pi}{(1+n_1)\Omega}.
\end{equation}
\par

The improvements obtained when using the higher order analytical solution, with short-period corrections up to the order eight, instead of the lower order analytical solution, with short-period corrections up to the sixth order, are not too relevant for orbits with small amplitude libration. However, one should note that the solution given by Eqs.~(\ref{qLin})--(\ref{QLin}) is, in fact, a ninth-order solution, because the ninth-order term of the averaged Ham\-il\-to\-ni\-an of the perturbation theory vanishes. Then, the higher order analytical solution can be extended to encompass the short-period corrections up to the order nine.
\par

\begin{figure}[htb]
\begin{center}
\includegraphics[scale=0.7, angle=0]{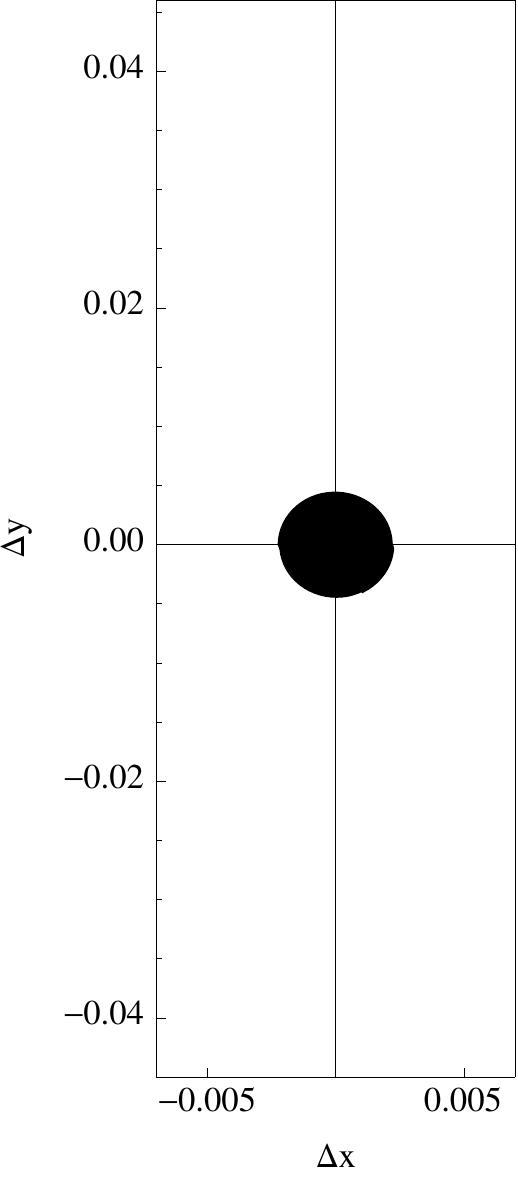}
\caption{Case 1 of Table \protect\ref{t:test6}. Errors of the 9th order analytical solution with short period corrections. Compare with Fig.~\ref{f:error1six}.}
\label{f:error1nine}
\end{center}
\end{figure}

Inclusion of the ninth order short-period corrections definitely improves the accuracy of the analytical solution. These improvements are clearly noticed when comparing the errors shown in Fig.~\ref{f:error1nine} for the test case 1 with those obtained with the low order analytical solution, shown in Fig.~\ref{f:error1six}. Now, the errors relative to the size of the orbit fall to the order of just a few tens of thousandths. 
\par

\begin{figure}[htb]
\begin{center}
\includegraphics[scale=0.7, angle=0]{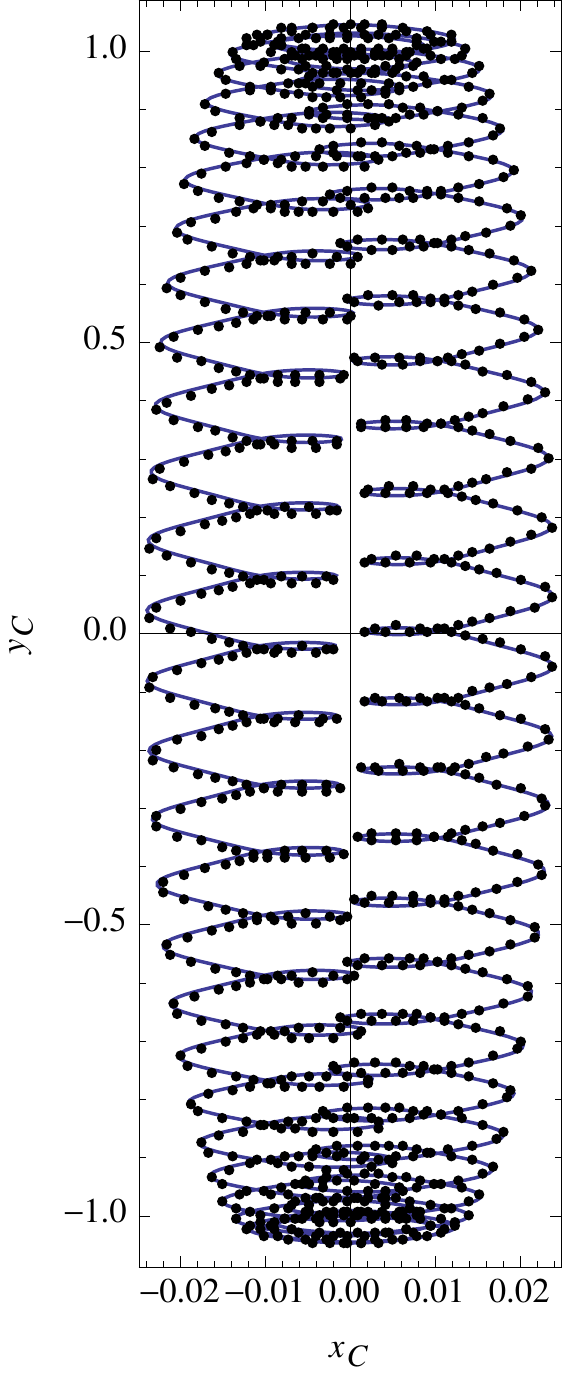}
\caption{Case 2 of Table \protect\ref{t:test6}. Trajectory $(x_C,y_C)$ of the center of the reference ellipse of the ninth-order analytical solution, with short-period corrections included, (black dots) superimposed to the true solution (full line) along one libration period. Compare with Fig.~\protect\ref{f:orbit2sixce}.}
\label{f:orbit2ninece}
\end{center}
\end{figure}

Analogous improvements are also found in test case 2. They are clearly noticed when comparing Fig.~\ref{f:orbit2ninece} with the right plot of Fig.~\ref{f:orbit2sixce}. Now, the shift between the center of the reference ellipse of the true orbit and the corresponding one computed from the analytical solution is very small, and both trajectories almost match.
\par

\begin{figure}[htb]
\begin{center}
\includegraphics[scale=0.8, angle=0]{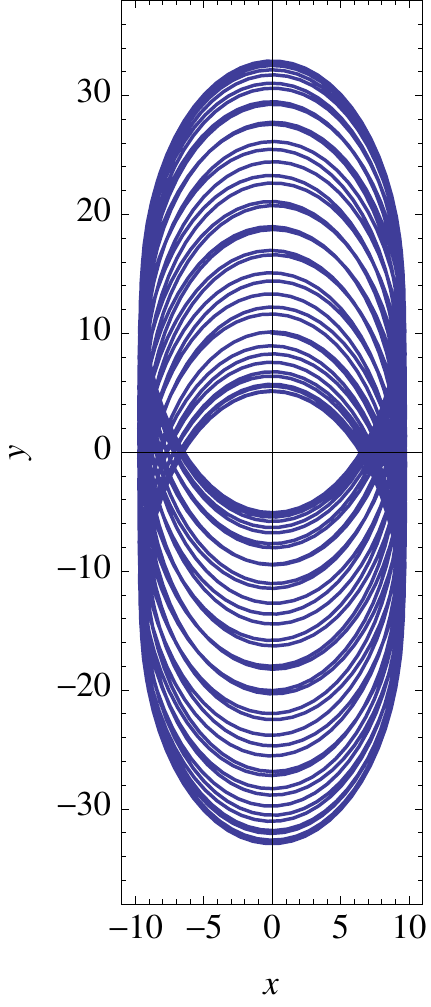} \qquad
\includegraphics[scale=0.8, angle=0]{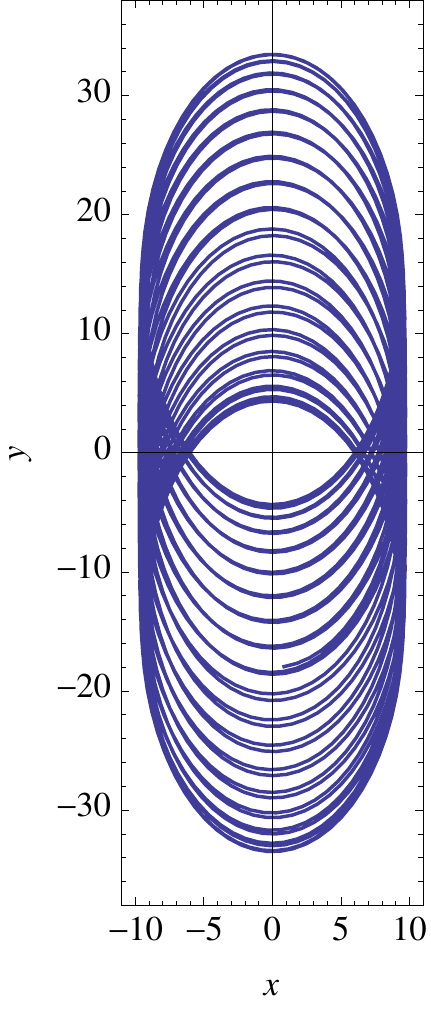} \qquad
\caption{Same as Fig.~\protect\ref{f:orbit3six}, but now the left plot is computed from the ninth-order analytical solution, which predicts a shorter libration period.}
\label{f:orbit3eigth}
\end{center}
\end{figure}

On the other hand, the 9th order analytical solution clearly improves the accuracy of orbits that get close to the primary, as illustrated in Fig.~\ref{f:orbit3eigth}, for the orbit of the test case 3 of Table \ref{t:test6}. Now, the amplitude of the librations predicted by the analytical solution are much closer to the actual ones, and the predicted orbit is much more similar in shape to the actual orbit. 
\par

\begin{figure}[htb]
\begin{center}
\includegraphics[scale=0.7, angle=0]{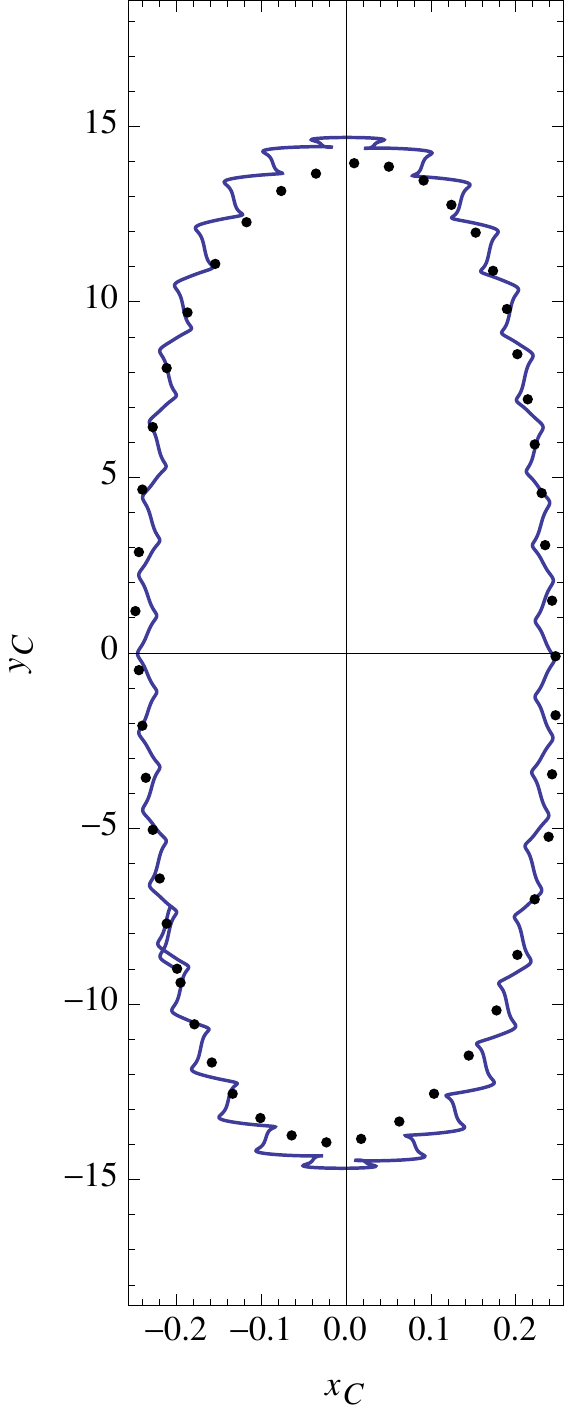}
\caption{Same as Fig.~\protect\ref{f:orbit3sixce} computed from the ninth order perturbation solution. Note the different scales of the $x_C$ and $y_C$ axes.}
\label{f:orbit3eigthce}
\end{center}
\end{figure}

\begin{figure}[htb]
\begin{center}
\includegraphics[scale=0.8, angle=0]{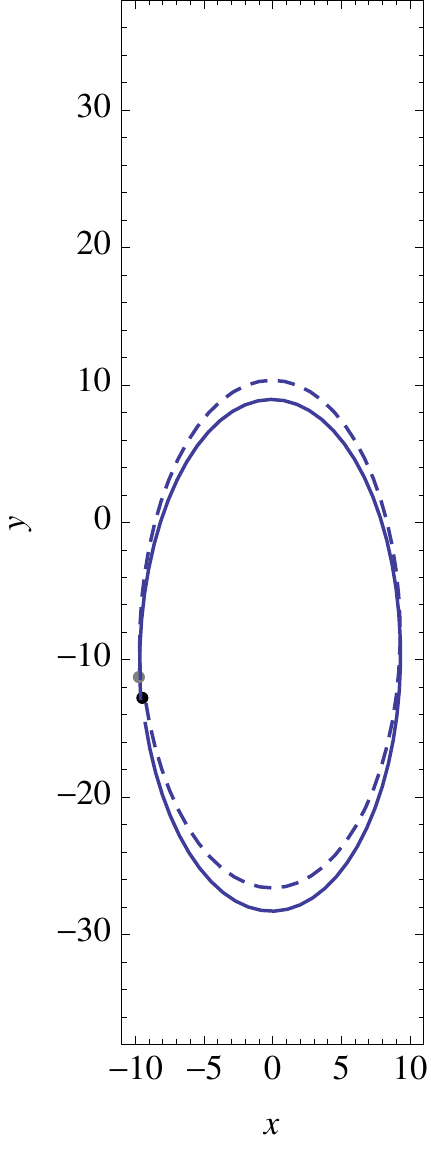}
\includegraphics[scale=0.8, angle=0]{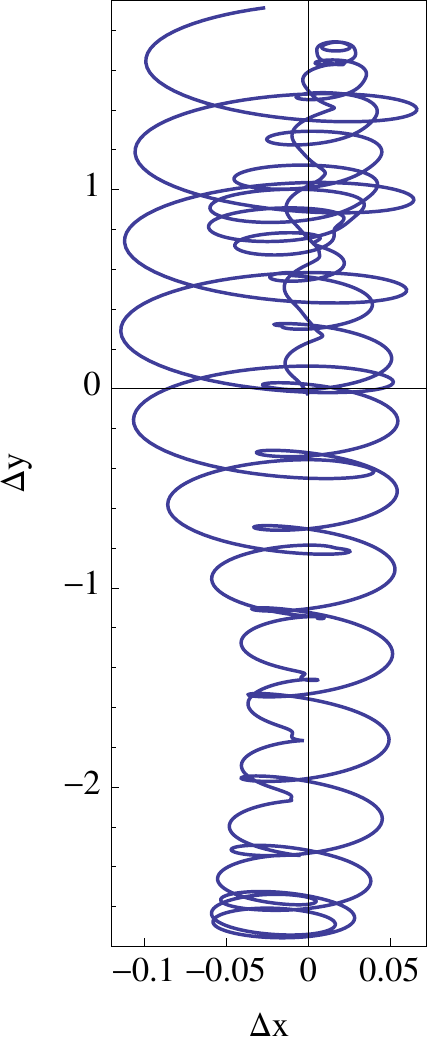}
\caption{Left: final part of the ninth order analytical solution (full line) and the true orbit (dashed line) in Fig.~\protect\ref{f:orbit3eigth}. Right: errors between both approaches along one libration period.}
\label{f:orbit3eigtherr}
\end{center}
\end{figure}

Still, the ninth-order truncation of the perturbation theory misses relevant terms of the smaller primary gravitation in the case of close encounters with the massless body, and the amplitude of the librations predicted by the ninth-order analytical solution are slightly different from the real ones, as clearly observed in the evolution of the center of the reference ellipse of the test case 3 orbit presented in Fig.~\ref{f:orbit3eigthce}. Because of this mismodeling, the libration period of the center of the reference ellipse predicted by  Eq.~(\ref{periodl8})  $T^*=236.66$ is slightly larger than the the libration period of the actual orbit, $T^*\approx232.21$. This difference between the predicted and the actual libration periods introduces a secular drift between both reference ellipses, as noticed in the left plot of Fig.~\ref{f:orbit3eigtherr}, which, for clarity, only displays the last orbital period of the propagation. Dots in this plot at one end of each curve belong to the same time of the propagation, which corresponds to the initial point of the last orbital period. The differences between the coordinates of the ninth order analytical solution and those of the true orbit are presented in the right plot of Fig.~\ref{f:orbit3eigtherr}, which shows how the errors concentrate mainly in the $y$ axis direction (note the different scales of the abscissas and ordinates in this plot).
\par

Computing higher orders of the perturbation theory may improve the solution. But, the eventual appearance of special functions in the computation of the generating function of the Lie transformation introduces serious difficulties in the computation of higher orders. This fact establishes a practical limit in the correct modeling of the gravitational effect of the smaller primary for orbits with large librations. In particular the analytical integration of the incomplete elliptic integrals of the first and second kinds must be tackled at the order 10 of the perturbation theory, a fact that notably complicates the computation of corresponding short-period corrections. In our case, this transformation has been computed only up to the ninth order. 
\par

As an alternative to the analytical solution, the numerical integration of Eqs.~(\ref{dq8})--(\ref{dQ8}) 
is very fast and efficient when compared to the numerical integration of Eqs.~(\ref{xp})--(\ref{Yp}) in intervals of one or several libration periods, because it is free of short-period effects. This is the approach we take in what follows, where this kind of numerically integrated solution is called a perturbation solution. Still, this alternative approach does not solve the problem of improving the perturbation model for large libration orbits, which requieres the computation of higher orders of the perturbation theory.
\par

The trouble in modeling perturbations efficiently is not unusual in astrodynamical problems, where perturbations of different nature may alternate in assuming the leading role along one orbital period. As an example of this situation we mention the two important cases of planetary orbits and highly elliptical orbits (see \cite{Lara2008,Lara2010,LaraPalacianYanguasCorral2010,ArmellinSanJuanLara2015}, for instance). In the present case, we found partial remedy to the bad modeling of the primary gravitation in including the averaged effects of the terms $K_{m,0}$, $M<{m}\le{N}$, all together as if they all where perturbations of the $M$-th order.
\par

\begin{figure}[htb]
\begin{center}
\includegraphics[scale=0.7, angle=0]{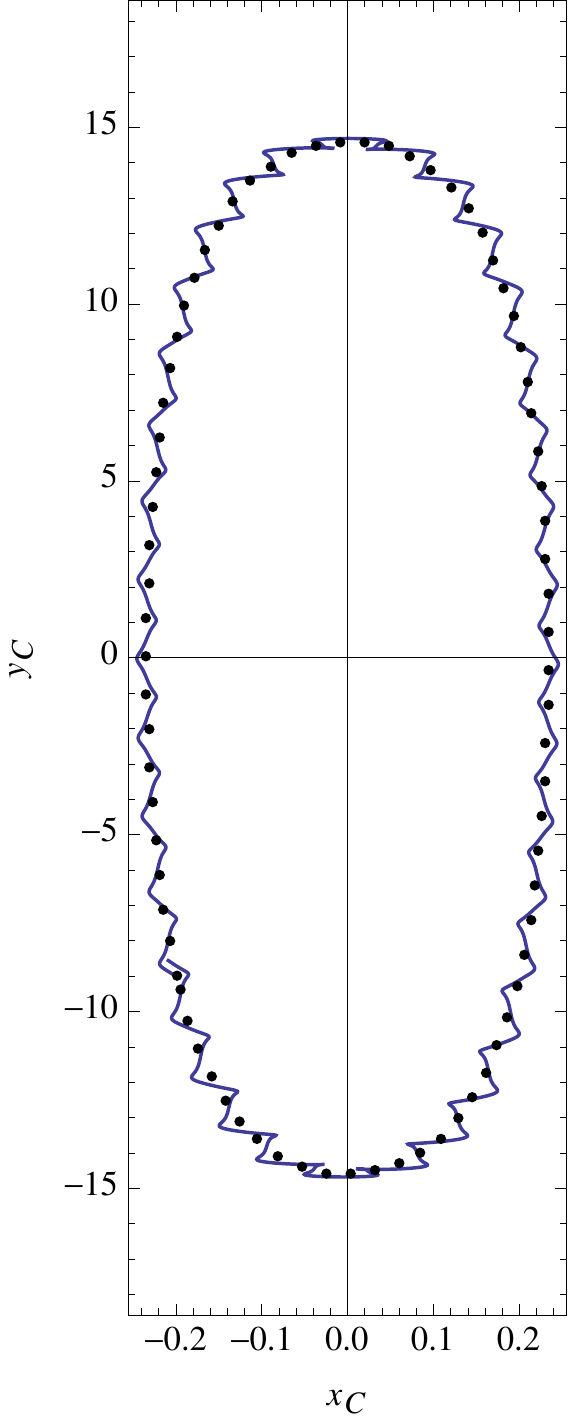} \qquad
\caption{Case 3 of Table \protect\ref{t:test6}. Trajectory of the center of the reference ellipse along the libration period $T^*=232.5$. The augmented 9th order perturbation solution (black dots) is superimposed to the true solution (full line).}
\label{f:orbit3eigthplus}
\end{center}
\end{figure}

\begin{figure}[htb]
\begin{center}
\includegraphics[scale=0.8, angle=0]{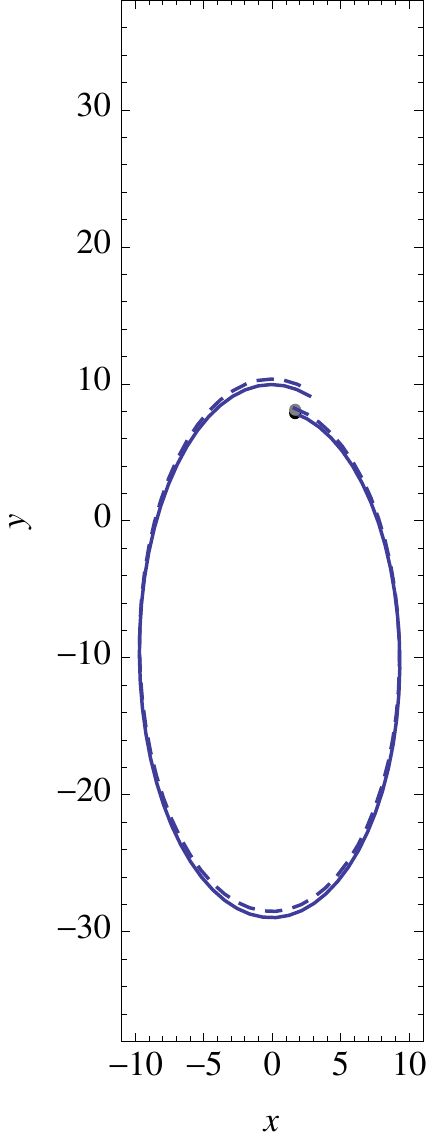}
\includegraphics[scale=0.8, angle=0]{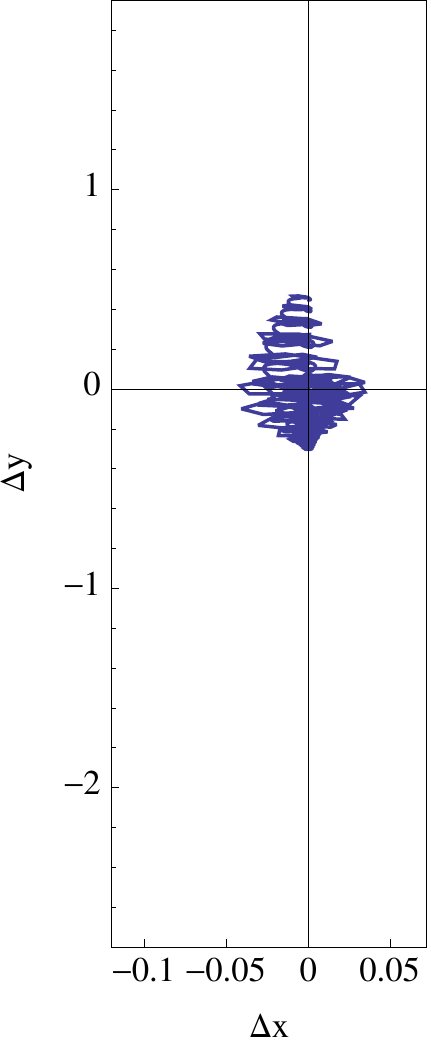}
\caption{Same as Fig.~\protect\ref{f:orbit3eigtherr} computed with the augmented ninth order perturbation solution.}
\label{f:err3eigthplus}
\end{center}
\end{figure}

The improvements obtained in the secular terms are illustrated in Figs.~\ref{f:orbit3eigthplus} and \ref{f:err3eigthplus} for the particular test case 3 of Table \ref{t:test6}. We set $M=9$ and $N=20$ ---which includes the averaged effects of $\sigma$ up to the order 8 and $\chi$ up to the order 16 in the expansion of the inverse of Eq.~(\ref{ro}). Now, the trajectory of the center of the reference ellipse predicted by the augmented 9th order perturbation solution clearly agrees, on average, with the numerical one, as shown in Fig.~\ref{f:orbit3eigthplus}. In consequence, the errors between the true solution and the perturbation one reduce substantially, as noted when comparing Fig.~\ref{f:err3eigthplus} and Fig.~\ref{f:orbit3eigtherr}. 
\par

\section{Conclusions}

Nonlinearities in the motion of quasi-satellite orbits with large librations are difficult to account for due to two main reasons. On the one hand, special functions which prevent closed form integration appear soon in a perturbation approach. On the other hand, the strength of the perturbations undergone by the massless body in its journey along the orbit is quite uneven due to the important variations of the distance to the smaller primary. 
Both difficulties are counterbalanced, yet only to some extent, with an efficient arrangement of the different orders of the perturbation of the Hill problem Hamiltonian.
\par

As a result of this Hamiltonian arrangement, a simple, low-order analytical perturbation solution shows that distant retrograde orbits are perturbed ellipses, whose centers evolve harmonically with a frequency that is proportional to the rotation rate of the system. This frequency gradually slows down for increasing values of the semi-minor axis of the reference ellipse. The low-order solution also shows that the linear growth of the phase of the massless body in the reference ellipse is modulated by long-period oscillations with half the period of the center of the reference ellipse's period. 
\par

Higher orders of the perturbation approach capture the non-linearities of the dynamics of orbits with large librations, but this is at the expense of loosing the simplicity of the analytical solution, which now requires the use of special functions and inversion procedures to be made explicit. Still, the oscillatory character of the librations is easily disclosed with simple contour plots of the averaged Hamiltonian, and analytical approximations can be computed by the Lindstedt-Poincar\'e method. Besides, the numerical integration of the averaged equations is very fast, and provides an efficient alternative to the purely analytical approach, which allows for increased accuracy in the propagation of orbits with large amplitude libration.

\subsection*{Acknowledgemnts}
Support by the Spanish State Research Agency and the European Regional Development Fund under Projects ESP2016-76585-R and ESP2013-41634-P (AEI/ERDF, EU) is recognized.

\appendix
\section{Short-period corrections} \label{a:LT}
The transformation equations of the Lie transforms procedure provide the short-period effects that were removed in the averaging. The direct transformation is given in the form
\[
\xi=\xi'+\sum_{m>0}\xi_{0,m}(\phi',q',\Phi',Q'),
\]
where $\xi$ stands for any of the variables $\phi,q,\Phi,Q$. For the Hamiltonian arrangement made, $\xi_{0,m}=0$ for $m<4$; the remaining coefficients up to $m=9$ are given below. We suppressed the prime notation without risk of confusion, and used the abbreviations $c\equiv\cos\phi$, $s\equiv\sin\phi$, and $\tilde{F}(\phi)=2K(3/4)\phi-F(\phi|3/4)$, $\tilde{E}(\phi)=2E(3/4)\phi-E(\phi|3/4)$, which are periodic functions of $\phi$ with period $2\pi$.
\small
\begin{eqnarray*}
\phi_{0,4} &=&-\frac{\mu}{\omega^2B^3}\frac{1}{2}\tilde{F}(\phi)
\\
q_{0,4} &=& 0 
\\
\Phi_{0,4} &=& \frac{\mu}{\omega{B}}\Big(\frac{1}{\Delta}-\tilde{K}\Big) 
\\
Q_{0,4} &=& 0
\\[1ex]
\phi_{0,5} &=& -\frac{\mu}{\omega^2B^3}\frac{2}{\Delta}\chi{s} 
\\
q_{0,5} &=& 0 
\\
\Phi_{0,5} &=& -\frac{\mu}{\omega{B}}\frac{4}{\Delta^3}\chi{c} 
\\
Q_{0,5} &=& -\frac{\mu}{\omega{B}^2}\frac{1}{2\Delta}{s}
\\[1ex]
\phi_{0,6} &=& \frac{\mu}{\omega^2B^3}\Big\{\frac{4}{\Delta}\sigma{c}
+\chi^2\Big[\tilde{E}(\phi)-\tilde{F}(\phi)
+\frac{3}{2\Delta^3}(5+3c^2)sc\Big]\Big\} 
\\
q_{0,6} &=& -\frac{\mu}{\omega^2B^2}\frac{2}{\Delta}{c}
\\
\Phi_{0,6} &=& \frac{\mu}{\omega{B}}\Big\{\chi^2\Big[\frac{2}{3}(\tilde{E}-\tilde{K})
-\frac{2}{\Delta^5}(1-9c^2)\Big]-\frac{2}{\Delta^3}\sigma{s}\Big\}
\\
Q_{0,6} &=& \frac{\mu}{\omega{B}^2}\chi\Big[\frac{1}{3}\big[\tilde{E}(\phi)-\tilde{F}(\phi)\big]
+\frac{1}{2\Delta^3}(5+3c^2)sc\Big]
\\[1ex]
\phi_{0,7} &=& \frac{\mu}{\omega^2B^3}\chi\Big[  8\sigma\Big(\frac{1}{\Delta^3}-\tilde{E}\Big)
+\frac{\chi^2}{3\Delta^5}(4-112c^2-84c^4){s} \Big]
\\
q_{0,7} &=& -\frac{\mu}{\omega^2B^2}\frac{8}{3}\chi\Big(\frac{1}{\Delta ^3}-\tilde{E}\Big) 
\\
\Phi_{0,7} &=& \frac{\mu}{\omega{B}}\frac{4}{\Delta^5}\chi
\Big[6\sigma{s}-\frac{\chi^2}{\Delta^2}(6-22c^2)\Big]c
\\
Q_{0,7} &=& \frac{\mu}{\omega{B}^2}\Big[ \frac{4}{3}\sigma\Big(\frac{1}{\Delta^3}-\tilde{E}\Big)
+\frac{\chi^2}{2\Delta^5}(1-28c^2-21c^4){s} \Big]
\\[1ex]
\phi_{0,8} &=& \frac{\mu}{\omega^2B^3}\Big\{\frac{\mu}{\omega^2B^3}
\Big[2\big[\phi+\tan^{-1}(2\cot\phi)\big]-\Big(\frac{3}{2}\tilde{K}+\frac{1}{2\Delta}\Big)\tilde{F}(\phi)\Big]
+\sigma^2\Big[\tilde{F}(\phi)-4\tilde{E}(\phi)-\frac{6}{\Delta^3}(2+3c^2)sc\Big] \\
&& -\frac{16}{\Delta^5}\sigma\chi^2 (3-5c^2-6c^4){c}
+\frac{5}{36}\chi^4\Big[14\tilde{E}(\phi)-11\tilde{F}(\phi)
-\frac{3}{\Delta^7}(67-387c^2-387c^4-189c^6)sc
\Big]\Big\}
\\
q_{0,8} &=& \frac{\mu}{\omega^2B^2}\Big\{
\frac{2}{3}\sigma\Big[4\tilde{E}(\phi)-\tilde{F}(\phi)+\frac{6}{\Delta^3}(2+3c^2)sc\Big] 
+\frac{4}{\Delta^5}\chi^2 (3-5c^2-6c^4){c}\Big\}
\\
\Phi_{0,8} &=& \frac{\mu}{\omega{B}}\bigg\{ \frac{\mu}{\omega^2B^3}\Big[ \frac{1}{2}-\frac{1}{\Delta^2}
+\tilde{K}\left(\frac{1}{\Delta}-\tilde{K}\right)-\frac{3}{2\Delta^3}\tilde{F}(\phi)sc\Big]
+\sigma^2\Big[\frac{2}{3}(\tilde{K}-4\tilde{E})+\frac{4}{\Delta^5}(1-3c^2)\Big] \\
&& +\chi^4\Big[\frac{1}{18}(14\tilde{E}-11\tilde{K})+\frac{2}{\Delta^9}(3-102c^2+227c^4)\Big]
-\frac{12}{\Delta^7}\sigma\chi^2(1-17c^2){s} \Big\}
\\
Q_{0,8} &=& \frac{\mu\chi}{\omega{B}^2}
\Big\{\chi^2\Big[\frac{14\tilde{E}(\phi)-11\tilde{F}(\phi)}{18}
+\frac{1}{6\Delta^7}(67-387c^2-387c^4-189c^6)sc\Big]
-\frac{4\sigma}{\Delta^5}(3-5c^2-6c^4){c} \Big\}
\\[1ex]
\phi_{0,9} &=& \frac{\mu}{B^3\omega^2}\Big\{\sqrt{3}\sigma\sinh^{-1}\big(\sqrt{3}{c}\big)
+ \frac{\mu}{\omega^2B^3}\chi\Big[ \frac{5}{8}\sqrt{3}\log\Big(\frac{2\sqrt{3}+3{s}}{2\sqrt{3}-3{s}}\Big)
+\frac{1}{\Delta}\Big(\frac{11}{2\Delta}-8\tilde{K}\Big){s} +\frac{4}{\Delta^3}\tilde{F}(\phi){c}\Big] \\
&& -\frac{2\chi}{\Delta^5}\Big[ \frac{40}{9\Delta^2}\sigma\chi^2(1-57c^2)
+8\sigma^2(1-4c^2-3c^4){s}
+\frac{\chi^4}{5\Delta^4}\big[19-634c^2+891(2+2c^2+c^4)c^4\big]s \Big] \Big\}
\\
q_{0,9} &=& -\frac{\mu}{\omega^2B^2}\Big\{ \frac{\sqrt{3}}{2}\sinh^{-1}\big(\sqrt{3}{c}\big) 
-\frac{8}{\Delta^5}\chi\Big[\frac{2\chi^2}{9\Delta^2}(1-57c^2)+
  \sigma(1-4c^2-3c^4){s} \Big] \Big\}
\\
\Phi_{0,9} &=& \frac{\mu}{\omega{B}}\frac{1}{\Delta^4}\Big\{ \frac{\mu}{\omega^2B^3}\chi\Big[
\Big(\frac{21}{2}-\frac{26}{\Delta}\tilde{K}\Big){c}
-\Big(\frac{3}{2}-\frac{6}{\Delta}\tilde{K}\Big)(3-4c^2){c}
-\frac{2}{\Delta}F(\phi)(1-6c^2){s} \Big]
-\frac{3}{2}\Delta^3\sigma{s} \\
&& -\frac{8}{\Delta^3}\chi\Big[12\sigma^2(1-2c^2)
+\frac{\chi^4}{\Delta^4}(15-190c^2+303c^4) 
+\frac{10}{\Delta^2}\sigma\chi^2(3-19c^2){s}\Big]{c} \Big\}
\\ 
Q_{0,9} &=& \frac{\mu}{\omega{B}^2}\Big\{ \frac{\mu}{\omega^2B^3}\Big[ 
\frac{\sqrt{3}}{16}\log\Big(\frac{2\sqrt{3}+3{s}}{2\sqrt{3}-3{s}}\Big)
+\frac{1}{\Delta}\Big(\frac{1}{\Delta^2}\tilde{F}(\phi){c}-\tilde{K}{s}\Big)
+\frac{3}{4\Delta^2}{s}\Big] \\
&& -\frac{1}{\Delta^5}\Big[\frac{8}{3\Delta^2}\sigma\chi^2(1-57c^2)
+2\sigma^2(1-4c^2-3c^4){s} 
+\frac{\chi^4}{6\Delta^4}\big[19-634c^2+891(2+2c^2+c^4)c^4\big]s \Big] \Big\}
\end{eqnarray*}
\par

The inverse transformation is given by
\[
\xi'=\xi+\sum_{m>0}\xi'_{0,m}(\phi,q,\Phi,Q),
\]
where $\xi'\in(\phi',q',\Phi',Q')$. Up to the order 9, terms $\xi'_{0,m}$ are formally the opposite of $\xi_{0,m}$ except for
\begin{eqnarray*}
\phi'_{0,8} &=& \frac{\mu^2}{\omega^4B^6}\left(\frac{1}{\Delta}-\tilde{K}\right)\tilde{F}(\phi)-\phi_{0,8}
\\
\Phi'_{0,8} &=& \frac{\mu^2}{\omega^3B^4}\Big[
\tilde{K}\Big(\frac{2}{\Delta}-\tilde{K}\Big)-\frac{1}{\Delta ^2}+\frac{3}{2\Delta^3}\tilde{F}(\phi)sc \Big]
-\Phi_{0,8}
\\[1ex]
\phi'_{0,9} &=& \frac{\mu^2}{\omega^4B^6}\frac{2\chi}{9\Delta}
\Big\{\frac{1}{\Delta^2}\tilde{F}(\phi){c}
+37\Big(\frac{1}{\Delta}-\tilde{K}\Big){s}\Big\}
-\phi_{0,9}
\\
\Phi'_{0,9} &=&\frac{\mu^2}{\omega^3B^4}\frac{\chi}{\Delta}\Big[
\frac{29}{9}\Big(\frac{1}{\Delta}-\tilde{K}\Big){c}
+\Big(\frac{5}{\Delta}-\frac{29}{3}\tilde{K}\Big)\frac{1}{\Delta^2}{c}\sin^2\phi
-\frac{22}{9\Delta^4}\tilde{F}(\phi)(1-6c^2){s} \Big] -\Phi_{0,9}
\\
Q'_{0,9} &=& \frac{\mu^2}{\omega^3B^5}\frac{11}{9\Delta}
\Big[\frac{1}{\Delta^2}\tilde{F}(\phi){c}+\Big(\frac{1}{\Delta}-\tilde{K}\Big){s} \Big] -Q_{0,9}
\end{eqnarray*}

\end{document}